\documentclass[12pt,leqno]{article}
\usepackage{upgreek}

\usepackage{xcolor}
\usepackage{bm,amsmath,amsthm,hyperref,amsfonts,amssymb,color}
\usepackage{mathrsfs} 
\usepackage{cite}

\usepackage{relsize} 
\usepackage{psfrag}  

\usepackage{tikz}  
\usetikzlibrary{decorations.pathmorphing}
\usetikzlibrary{patterns}
\usepackage{xcolor} 
\usepackage{pgfplots}
\usetikzlibrary{calc}

\newcommand{\EPS}{\varepsilon}

\newcommand\DT[1]{\mathchoice
                 {{\buildrel{\hspace*{.1em}\text{\LARGE.}}\over{#1}}}
                 {{\buildrel{\hspace*{.1em}\text{\LARGE.}}\over{#1}}}
                 {{\buildrel{\hspace*{.1em}\text{\Large.}}\over{#1}}}
                 {{\buildrel{\hspace*{.1em}\text{\large.}}\over{#1}}}}
\newcommand\DTk[1]{\mathchoice
                 {{\buildrel{\hspace*{.1em}\text{\,\LARGE.}\EPS k\!\!\!\!\!}\over{#1}}}
                 {{\buildrel{\hspace*{.1em}\text{\,\LARGE.}\EPS k\!\!}\over{#1}}}
                 {{\buildrel{\hspace*{.1em}\text{\,\Large.}\EPS k\!\!}\over{#1}}}
                 {{\buildrel{\hspace*{.1em}\text{\,\large.}\EPS k\!\!}\over{#1}}}}
\newcommand\pdt[1]{\frac{\partial{#1}}{\partial t}} 
\newcommand{\lineunder}[2]{\LU{\begin{array}[t]{c}\underbrace{#1}\vspace*{.5em}\end{array}}{\mbox{\footnotesize\rm #2}}}
\newcommand{\LU}[2]{\begin{array}[t]{c}#1\vspace*{-1em}\\_{#2}\end{array}}
\newcommand{\linesunder}[3]{\LSU{\begin{array}[t]{c}\underbrace{#1}\vspace*{.5em}\end{array}}{\mbox{\footnotesize\rm #2}}{\mbox{\footnotesize\rm #3}}}
\newcommand{\LSU}[3]{\begin{array}[t]{c}#1\vspace*{-1em}\\_{#2}\vspace*{-.5em}\\_{#3}\end{array}}
\newcommand{\morelinesunder}[4]{\LSUU{\begin{array}[t]{c}\underbrace{#1}\vspace*{.5em}\end{array}}{\mbox{\footnotesize\rm #2}}{\mbox{\footnotesize\rm #3}}{\mbox{\footnotesize\rm #4}}}
\newcommand{\LSUU}[4]{\begin{array}[t]{c}#1\vspace*{-1em}\\_{#2}\vspace*{-.5em}\\_{#3}\vspace*{-.5em}\\_{#4}\end{array}}
\newcommand{\Item}[2]{\parbox[t]{.05\textwidth}{#1}\hfill%
      \parbox[t]{.95\textwidth}{#2}\vspace*{.8mm}} 
\newcommand{\divS}{\mathrm{div}_{\scriptscriptstyle\textrm{\hspace*{-.1em}S}}^{}}
\newcommand{\nablaS}{\nabla_{\scriptscriptstyle\textrm{\hspace*{-.3em}S}}^{}}

\def\Vdots{\!\mbox{\setlength{\unitlength}{1em}
\begin{picture}(0,0)
\put(-.07,0){.}
\put(-.07,.3){.}
\put(-.07,.6){.}
\end{picture}\hspace*{.2em}}}
\usepackage{mathrsfs}   
\usepackage{eucal}      

  \def\bbI{{\mathbb I}}

\def\FG{\boldsymbol}
   
 \def\ee{{\FG e}}  
  
\def\jj{{\FG j}}   
 \def\nn{{\FG n}}

\def\vv{{\FG v}}  \def\xx{{\FG x}}

\def\DD{{\FG D}} 
\def\FF{{\FG F}}

 \def\TT{{\FG T}} 

\newcommand{\R}{\mathbb R}
\newcommand{\N}{\mathbb N}
\newcommand{\Nabla}{{\nabla}}
\newcommand{\Fe}{\FF}
\newcommand{\FFeps}{\FF_{\!\EPS}}

\newcommand{\FFepsk}{\FF_{\!\EPS k}}

\newcommand{\eq}[1]{(\ref{#1})}
\renewcommand{\d}{\mathrm d}  
\newcommand{\barOmega}{\hspace*{.2em}{\overline{\hspace*{-.2em}\varOmega}}}

\newtheorem{theorem}{Theorem}[section]

\newtheorem{definition}[theorem]{Definition}
\newtheorem{example}[theorem]{Example}

\newtheorem{remark}[theorem]{Remark}

\numberwithin{equation}{section}

\def\upyy{\text{\bf{y}}}
\def\upxx{\text{\bf{X}}}
\def\upvv{\text{\bf{v}}}
\def\upFF{\text{\bf{F}}}

\begin{document}
\begin{sloppypar}

\allowdisplaybreaks

\begin{center}{\Large\bf 
Thermo-elastodynamics of finitely-strained multipolar
\\[.3em]viscous solids with an energy-controlled stress.
}

\bigskip\bigskip

{\sc Tom\'a\v s Roub\'\i\v cek}
\footnote{Support from CSF grant no.\,23-06220S, and from the institutional
support RVO:61388998 (\v CR)  is acknowledged.}

\end{center}

\bigskip
\noindent
{Mathematical Institute, Charles University,\\
Sokolovsk\'a 83, CZ-186~75~Praha~8, Czech Republic,
\\[.2em]
and\\[.2em]
Institute of Thermomechanics of the Czech Academy of Sciences,\\ 
Dolej\v skova 5, CZ-182 00 Praha 8, Czech Republic},\\[.2em]
email: {\tt tomas.roubicek@mff.cuni.cz}, {\smaller ORCID 0000-0002-0651-5959}.

\date{}

\begin{abstract}
The thermodynamical model of viscoelastic deformable solids
at finite strains with Kelvin-Voigt rheology with a higher-order viscosity
(using the concept of multipolar materials) is formulated in
a fully Eulerian way in rates. Assumptions used in this paper
allow for a physically justified free energy leading to
non-negative entropy that satisfies the 3rd law of thermodynamics,
i.e.\ entropy vanishes at zero temperature, and energy-controlled stress.
This last attribute is used advantageously to prove the existence
and a certain regularity of weak solutions by a simplified
Faedo-Galerkin semi-discretization, based on estimates obtained
from the total-energy and the mechanical-energy balances.
Some examples that model neo-Hookean-type materials are presented, too.

  \medskip


\medskip

\end{abstract}



\def\TRACTION{\bm{f}}
\def\GRAVITY{\bm{g}}
\def\rhoR{\varrho_\text{\sc r}^{}}
\def\RRR{\text{\sc r}}
\def\LAM{\lambda}
\def\W{w}
\def\OMEGA{\omega}
\def\ALPHA{\alpha}
\def\ZETA{\gamma}
\def\ZETA{a}
\def\ONEALPH{2}
\def\DIS{\DD}
\def\TWO{2}
\def\TWOprime{2}
\def\EXP{\mu}
\def\wh{\widehat}
\def\wt{\widetilde}

\def\Eng{e}
\def\ENG{E}
\def\Ent{u}
\def\ENT{U}
\def\Sv{{\bm D}}
\def\Colon{{:}}

\section{Introduction}

The heat-transfer problem is an ``evergreen'' problem  studied  in the
mathematical literature for centuries, in recent decades in 
non-linear variants in a mechanical context and possibly coupled
with other physical or chemical processes. Often (or even
mostly) the heat capacity is considered to be ``uniformly'' positive
(quite typically just constant), which is only possible if the entropy
decays to $-\infty$ as the temperature approaches absolute zero.
Although the constant heat capacity is the first choice in most
mathematical and engineering publications,
such a choice is not physical because the entropy $-\infty$ at zero
temperature contradicts the sound physical arguments. Rather, the heat
capacity should degenerate to zero at zero temperature, so that
the entropy remains finite (usually calibrated to zero)
at zero temperature, which is a physically desirable
(and even ultimate) property.

This degeneracy need not be detrimental to the mathematical analysis.
In this article, we want to revisit the {\it heat-transfer problem at
finite strains} in the {\it Eulerian framework}, which is probably the most
complicated due to $L^1$-heat sources arising from the mechanical
visco-elastodynamic compressible model together with transport equations
and convective terms in the heat equation. On the other hand, we confine
ourselves on the mathematically simplest viscoelastic {\it rheology} of
the {\it Kelvin-Voigt} type.

The Eulerian formulation  itself  is relatively standard,
cf.\ \cite{GuFrAn10MTC,Mart19PCM,Rubi21CMEF,Silh97MTCM,True69RT}.
Mathematical analysis in the conventional simple-material variant
is, however, open (as also articulated in \cite{Ball84MELE,Ball02SOPE})
except in the fluidic variant where the absence of any shear elastic response
 has been  exploited, cf.\ \cite{FeiNov09SLTV}. For completeness, we report
some existing analytical results in isothermal incompressible situations in
\cite{LeLiZh08GSIV,LiuWal01EDFC} and compressible in
\cite{HuMas16GSRH,QiaZha10GWPC}  with small data or with (non-physical)
convex stored energies. There seems to be a certain agreement that some
higher gradients (falling into the concept of so-called
{\it non-simple materials}) are desirable or even necessary for a reasonable
analysis (unless some very week e.g.\ of measure-valued-type solution concepts
like e.g.\ in Lagrangian formulation in
\cite{ChGaTz19MVSE,Demo00WSCN,DeStTz01VAST,GaViKo22WSUM} are used).
It should be openly admitted that  these analytical objectives are  the main
purpose  of including such higher gradients  in this paper, although
 a certain additional motivation could be
to open up options for modelling the dispersion of possible elastic waves
in some desirable way. In the Eulerian frame, such higher gradients are
to be involved rather in the dissipative than the conservative part (giving
rise continua which are sometimes referred as {\it multipolar}), so
that their influence is manifested itself only in fast evolution. Then the
higher gradients can lead to easier propagation of elastic waves with less
dispersion and less attenuation (compared with the usual first-order gradients)
especially at frequencies close to the critical frequency above which the
waves cannot propagate at all, as pointed out in \cite[Sect.3.1]{Roub23GTLV}
in the 1-dimensional linear variant.  Such higher gradients have been used in
\cite{Roub24TVSE} for essentially the same thermomechanical system as in this paper,
where a similar analysis by a semi-Galerkin approximation was used but under
rather strong and not entirely physical assumptions with limited applications
and exploiting a more complicated regularization.  

The main features  of the present article  can be summarized as follows:
\begin{itemize}
\vspace*{-.2em}\item
usage of the {\it actual free energy}, so that formally (as opposed to the referential
energy as in Remark~\ref{rem-engineering}) one does not need to control
explicitly the determinant of the deformation gradient (although the control
of the mass density implicitly includes this determinant),
\vspace*{-.8em}\item
estimation based on the {\it energy-controlled stress}, i.e.\ the Cauchy stress
whose magnitude is dominated by the actual internal energy,
\vspace*{-.8em}\item
sufficiently weak qualification of data which, in particular, admits physically
justified free energies
\begin{itemize}
\vspace*{-.8em}\item[$\circ$] satisfying the 3rd law of thermodynamics, i.e.\
yielding the non-negative entropy which vanishes at zero temperature (and is
then independent of the mechanical state), although this attribute will not be
explicitly exploited in the analysis, and
\vspace*{-.4em}\item[$\circ$] allowing for a natural continuous extension for
negative values of temperature (as in Fig.\,\ref{fig-free-energy} 
below) exploited advantageously for the Galerkin approximation.
\end{itemize}
\end{itemize}
In comparison with \cite{Roub24TVSE}, the mentioned weaker data qualification
will  allow  free energies describing e.g.\ thermal expansion and complying
with a heat capacity which remains positive for arbitrarily large deformations, not
only for conditionally not too large  ``moderate''  deformations.
Moreover, the semi-Galerkin approximation scheme used here is much simpler than in
\cite{Roub24TVSE} and in particular does not rely on any cut-off and the
non-negativity of the limit temperature is obtained more easily, among various
other simplifications.

The plan for this paper is as follows: In Section~\ref{sec-system},
we formulate the thermo-visco-elastic system at finite strains, and
briefly derive its thermomechanical consistency, i.e.\ 1st and 2nd laws
of thermodynamics. Then, in Section~\ref{sec-anal}, we present its
rigorous analysis as far as the global existence and certain regularity
of weak solutions for large data. Finally, in Section~\ref{sec-exa},
we present some examples of neo-Hookean materials complying with the
assumptions imposed in Section~\ref{sec-anal}.

Let us emphasize that the approximation and estimation strategy developed
in Section~\ref{sec-anal} seems to be well applicable to the (in some sense simpler)
variant in the Lagrangian framework and, a-fortiori, to a small-strain variant.

\section{The thermomechanical system and its energetics}\label{sec-system}

We briefly remind the fundamental concepts and formulas that can
mostly be found in various monographs, as e.g.\
\cite{GuFrAn10MTC,Mart19PCM,Rubi21CMEF,Silh97MTCM,True69RT}.

In finite-strain continuum mechanics, the basic geometrical concept is the
time-evolving deformation $\upyy:\varOmega\to\R^d$ as a mapping from a reference
configuration of the body $\varOmega\subset\R^d$ into a physical space $\R^d$.
The ``Lagrangian'' space variable in the reference configuration is
denoted as $\upxx\in\varOmega$ while in the ``Eulerian'' physical-space
variable as $\xx\in\R^d$. The basic geometrical object is the Lagrangian
deformation gradient $\upFF=\Nabla_{\!\upxx}^{}\upyy$.
We will be interested in deformations $\xx=\upyy(t,\upxx)$ evolving in time,
which are sometimes called ``motions''. The referential velocity
is $\upvv=\pdt{}\upyy$. The inverse mapping $\bm\xi:\xx\mapsto\upxx=\upyy^{-1}(t,\xx)$
is called a {\it return} (sometimes called also a {\it reference})
{\it mapping}.

Having $\bm\xi$, we define the Eulerian velocity $\vv=\upvv{\circ}\bm\xi$
and the Eulerian deformation gradient $\FF=\upFF{\circ}\bm\xi$. Having the
Eulerian velocity $\vv$, we define the {\it convective time derivative}
applied to scalars or, component-wise, to vectors or tensors,
using the dot-notation $(\cdot)\!\DT{^{}}$. By the chain-rule calculus,
we have the {\it transport equation-and-evolution for the
deformation gradient} and its determinant and its inverse as
\begin{align}\nonumber\\[-2.7em]
\DT\FF=(\nabla\vv)\FF\,, \ \ \ \ \ \ \ \
\DT{\overline{\det\FF}}=(\det\FF){\rm div}\,\vv\,,
\ \ \ \text{ and }\ \ \ \
\DT{\overline{\!\!\!\bigg(\frac1{\det\FF}\bigg)\!\!\!}}\ 
=-\frac{{\rm div}\,\vv}{\det\FF}\,.
\label{ultimate}\end{align}

The reference mapping $\bm\xi$, which itself is well defined by its transport
equation $\DT{\bm\xi}=\bm0$, actually does not explicitly occur in the
formulation of the problem if the medium is considered homogeneous in its
reference configuration. Here, we will benefit from the boundary condition
$\vv{\cdot}\nn=0$ below, which causes that the actual domain $\varOmega$ does
not evolve in time. The same convention applies to temperature $\theta$ and
thus also  the (conservative) Cauchy
stress  $\TT$,  the enthalpy  $\eta$, and  the dissipative
stress  $\DIS$ in \eq{stress-entropy} and \eq{Cauchy-dissip} below, which
will make the problem indeed fully Eulerian. Cf.\ the continuum-mechanics
textbooks as e.g.\
\cite{GuFrAn10MTC,Mart19PCM}. The mass density (in kg/m$^3$) transport and the
``mass sparsity'' as the inverse mass density $1/\varrho$ transport write: 
\begin{align}\nonumber\\[-2.7em]
\DT\varrho=-\varrho\,{\rm div}\,\vv\ \ \ \ \text{ and }\ \ \ \ 
\DT{\overline{\!\!\bigg(\frac1\varrho\bigg)\!\!}}\
=\frac{{\rm div}\,\vv}{\varrho}\,.
\label{cont-eq+}\end{align}
The former equation in \eq{cont-eq+} called the {\it continuity equation}
equivalently writes $\pdt{}\varrho+{\rm div}(\varrho\vv)=0$ (expressing
that the conservation of mass) and  ensures the transport of the momentum
$\varrho\vv$: 
\begin{align}\label{inertial}
\frac{\partial}{\partial t} (\varrho \bm{v})+
\text{\rm div}(\varrho \bm{v}{\otimes} \bm{v})=\varrho\DT\vv\,.
\end{align}
One can determine the density $\varrho$ instead of the transport equation for mass
density \eq{cont-eq+} from the algebraic relation
\begin{align}\nonumber\\[-2.7em]
\varrho=\frac{\rhoR}{\det\FF}\,,
\label{density-algebraically}\end{align}
where $\rhoR$ is the mass density in the reference
configuration. Later, we will consider the initial conditions $\FF_0$ for
\eq{ultimate} and $\varrho_0$ for \eq{cont-eq+}. 
These two should be related by $\varrho(0)=\rhoR/\!\det\FF_0$.

The main ingredients of the model are the (volumetric) {\it free energy}
$\psi$ depending on deformation gradient $\FF$ and temperature $\theta$,
and the temperature-dependent {\it dissipative stress} $\DIS$ not necessarily
possessing any underlying potential. The free energy $\psi=\psi(\FF,\theta)$
considered per the {\it actual volume} (in contrast to a referential energy as
in Remark~\ref{rem-engineering} below) leads to the conservative part of the
(actual) Cauchy stress, the entropy, and the  heat capacity respectively
\begin{align}
\TT=\psi_\Fe'(\Fe,\theta)\Fe^\top\!+\psi(\Fe,\theta)\bbI\,,\ \ \ \ \ \ 
\eta=-\psi_\theta'(\Fe,\theta)\,,\ \ \text{ and }\ \  
c(\FF,\theta)=-\theta\psi_{\theta\theta}''(\FF,\theta)\,,
\label{stress-entropy}\end{align}
where $\bbI$ denotes the unit matrix. The concept of stresses governed by
some potential is called {\it hyperelasticity}.

In the already anticipated visco-elastodynamic {\it Kelvin-Voigt rheological
model}, we will also use a dissipative contribution to the Cauchy stress,
which will make the system of a parabolic type compared to mere elastodynamics
which would be of a hyperbolic type. In addition to the usual first-order stress
depending on $\ee(\vv)$, we consider a dissipative contribution to the Cauchy
stress involving also a higher-order 2nd-grade {\it hyper-stress} $\mathcal{H}$,
so that the overall dissipative stress is:
\begin{align}
\DD=\nu_1\ee(\vv)-{\rm div}\,\mathcal{H}
\ \ \ &\text{ with }\ \ \text{ and }\ 
\mathcal{H}=\nu_2|\nabla^2\vv|^{p-2}\nabla^2\vv\,.
\label{Cauchy-dissip}\end{align}
Rather for notational simplicity, we consider $\nu_1>0$ and $\nu_2>0$ constant,
even though it would not pose any analytical difficulties to consider them
continuously dependent on $(\FF,\theta)$ and/or tensor-valued. The mentioned concept
of the 2nd-grade ``hyper-stress'' refers to a 3rd-order tensor (here denoted by
$\mathcal{H}$) whose divergence yields a contribution to the the Cauchy stress.
This modelling concept is inspired by R.A.\,Toupin \cite{Toup62EMCS} and
R.D.\,Mindlin \cite{Mind64MSLE} and has quite widely been used for the general
nonlinear context of {\it multipolar materials}, both fluids by J.~Ne\v cas at al.\
\cite{Neca94TMF,NeNoSi91GSCI,NecRuz92GSIV} or solids \cite{Ruzi92MPTM,Silh92MVMS}.

The {\it momentum equilibrium} equation then balances the divergence of
the total Cauchy stress with the inertial and gravity force:
\begin{align}
\varrho\DT\vv-{\rm div}\big(\TT{+}\DD)=\varrho\GRAVITY
\label{Euler-thermodynam1-}\end{align}
with $\TT$ from \eq{stress-entropy} and $\DD$ from \eq{Cauchy-dissip}. 

The second ingredient in \eq{stress-entropy} is subjected to 
the {\it entropy equation}:
\begin{align}
  \pdt\eta+{\rm div}\big(\vv\,\eta\big)
  =\frac{\xi-{\rm div}\,{\bm j}}\theta\ \ \ \ \text{ with }\ 
  \jj=-\kappa(\Fe,\theta)\nabla\theta \,
\label{entropy-eq}\end{align}
with $\xi$ the heat production rate  specified later.
The latter
equality in \eq{entropy-eq} is the {\it Fourier law} determining phenomenologically
the heat flux $\jj$ proportional to the temperature gradient through the thermal
conductivity coefficient $\kappa=\kappa(\Fe,\theta)$. Assuming $\xi\ge0$ and
$\kappa\ge0$ and integrating \eq{entropy-eq} over the domain $\varOmega$
while imposing the impenetrability of the boundary in the sense that
the normal velocity $\vv{\cdot}\nn$ vanishes across the boundary
$\varGamma$ of $\varOmega$, we obtain the {\it Clausius-Duhem inequality}:
\begin{align}\label{entropy-ineq}
\frac{\d}{\d t}\!\!\!\!\!\!\lineunder{\int_\varOmega\eta\,\d\xx}{total entropy}\!\!\!\!\!\!
&=\int_\varOmega\!\!\!\!\!\!\lineunder{\frac\xi\theta+\kappa(\FF,\theta)\frac{|\nabla\theta|^2\!}{\theta^2}}{entropy production rate}\!\!\!\!\!\d\xx
\\\nonumber
&+\int_\varGamma\!\!\!\!\lineunder{\Big(\kappa(\FF,\theta)\frac{\nabla\theta}{\theta}-\eta\vv\Big)}{entropy flux}\!\!\!\!\!\!\cdot\nn\,\d S
\ge\int_\varGamma\!\kappa(\FF,\theta)\frac{\nabla\theta\!}{\theta}{\cdot}\nn\,\d S\,.
\end{align}
If the system is thermally isolated in the sense that the normal heat flux
$\jj{\cdot}\nn$ vanishes across the boundary $\varGamma$, we recover the
{\it 2nd law of thermodynamics}, i.e.\ the total entropy in isolated
systems is nondecreasing in time.

The {\it internal energy} is given by the {\it Gibbs relation}
$\Eng=\psi+\theta\eta$. Using the calculus
\begin{align}\nonumber
 & \pdt\Eng+{\rm div}\big(\vv\,\Eng\big)=
 \DT\Eng+\Eng\,{\rm div}\,\vv
 =\ \DT{\overline{\psi(\FF,\theta)-\theta\psi_\theta'(\FF,\theta)}}
\,+\big(\psi(\FF,\theta)-\theta\psi_\theta'(\FF,\theta)\big)\,{\rm div}\,\vv
\\[-.2em]&\nonumber\ =
\psi'_{\Fe}(\FF,\theta){:}\DT\FF+\psi_\theta'(\FF,\theta)\DT\theta
-\theta\psi_{\Fe\theta}''(\FF,\theta)\DT\FF-\theta\psi_{\theta\theta}''(\FF,\theta)\DT\theta
-\DT\theta\psi_\theta'(\FF,\theta)
\\[-.2em]&\nonumber\hspace{22.3em}
+\big(\psi(\FF,\theta)-\theta\psi_\theta'(\FF,\theta)\big)\,{\rm div}\,\vv
\\[-.6em]&\nonumber\ =\psi'_{\Fe}(\FF,\theta){:}\DT\FF+\theta\DT{\overline{\eta(\FF,\theta)}}
+\big(\psi(\FF,\theta)-\theta\psi_\theta'(\FF,\theta)\big)\,{\rm div}\,\vv
\\[-.2em]&\nonumber\stackrel{\scriptsize\eq{entropy-eq}}{=}\!\psi'_{\Fe}(\FF,\theta){:}\DT\FF
+\xi-{\rm div}\,\jj-\theta\eta(\FF,\theta){\rm div}\,\vv
+\big(\psi(\FF,\theta)-\theta\psi_\theta'(\FF,\theta)\big)\,{\rm div}\,\vv
\\[-.2em]&\nonumber
\stackrel{\scriptsize\eq{ultimate}}{=}\!\psi'_{\Fe}(\FF,\theta)\FF^\top{:}\nabla\vv
+\xi-{\rm div}\,\jj+\psi(\FF,\theta){\rm div}\,\vv
\!\stackrel{\scriptsize\eq{stress-entropy}}{=}\!
\xi-{\rm div}\,\jj+\TT{:}\nabla\vv\,,
\end{align}
we obtain the {\it internal-energy equation}
\begin{align}\label{Euler-thermodynam3-}
\pdt\Eng+{\rm div}\big(\vv\Eng{+}{\bm j}\big)=\xi
+\TT{:}\nabla\vv\ \ &\text{ with }\ \xi=\nu_1|\ee(\vv)|^2\!+\nu_2|\nabla^2\vv|^p
\ \text{ and}
\\[-.2em]\ &\text{ with }\
\Eng=\ENG(\FF,\theta):=\psi(\Fe,\theta)-\theta\psi_\theta'(\Fe,\theta)
\,;
\nonumber\end{align}
{}the heat-production rate $\xi$ is here specified to fit the
total-energy balance \eq{Euler-engr-finite} later.

Altogether, let us formulate the thermo-visco-elastodynamic system for the
quadruple $(\varrho,\vv,\FF,\theta)$, composed from the equations \eq{ultimate},
\eq{cont-eq+}, \eq{Euler-thermodynam1-}, and \eq{Euler-thermodynam3-}:
\begin{subequations}\label{Euler-thermo-finite}\begin{align}
&\DT\varrho=-({\rm div}\,\vv)\varrho\,,
\label{Euler0-thermo-finite}
\\&\label{Euler1-thermo-finite}
\varrho\DT\vv={\rm div}\big(\TT(\FF,\theta){+}\Sv\big)
+\varrho\GRAVITY\ \ \ \
\text{ where }\ \ \TT(\FF,\theta)=\psi_\FF'(\FF,\theta)\FF^\top\!\!+\psi(\FF,\theta)\bbI
    \\
    &\hspace*{8em}
\text{ and }\ \ 
 \Sv=\nu_1
 \ee(\vv)-{\rm div}\,\mathcal{H}
 \ \ \text{ with }\ \ \mathcal{H}=\nu_2
 |\nabla^2\vv|^{p-2}\nabla^2\vv\,,
\nonumber
\\\label{Euler2-thermo-finite}
&\DT\FF=(\nabla\vv)\FF\,,
\\
&\label{Euler3-thermo-finite}
\DT\Eng={\rm div}\big(\kappa(\FF,\theta)\nabla\theta\big)
+\nu_1|\ee(\vv)|^2+\nu_2|\Nabla^2\vv|^p-({\rm div}\,\vv)\Eng
+\TT(\FF,\theta){:}\nabla\vv
\\&\hspace{8em}
\text{ with }\ \ \Eng=\ENG(\FF,\theta)\,,
\ \text{ where }\ \ENG(\FF,\theta):=
\psi(\FF,\theta)-\theta\psi_\theta'(\FF,\theta)\,.
\nonumber\end{align}\end{subequations}

Denoting by $\nn$ the unit outward normal to the (fixed) boundary $\varGamma$
of the domain $\varOmega$, we complete this system by suitable boundary
conditions: 
\begin{align}\label{Euler-thermodynam-BC}
&\vv{\cdot}\nn=0,\ \ 
\big[(\TT(\FF,\theta){+}\DD)\nn{-}\divS\big(\mathcal{H}
\nn\big)\big]_\text{\sc t}^{}\!\!
={\bm0},\ \ 
\Nabla^2\vv{:}(\nn{\otimes}\nn)={\bm0},\,\text{ and }\, 
\kappa(\FF,\theta)\nabla\theta{\cdot}\nn=h\!\!\!
\end{align}
with $[\,\cdot\,]_\text{\sc t}^{}$ a tangential part of a vector.
Here $\divS={\rm tr}(\nablaS)$ denotes the $(d{-}1)$-dimensional
surface divergence with ${\rm tr}(\cdot)$ 
the trace of a $(d{-}1){\times}(d{-}1)$-matrix and $\nablaS z=\nabla z-
(\nabla z{\cdot}\nn)\nn$ being the surface gradient of a field $z$. 
The first condition (i.e.\ normal velocity zero) expresses impenetrability of
the boundary, which  has already been used for \eq{entropy-ineq}
and  which is the  most frequently adopted  assumption  in
literature for the Eulerian formulation. This simplifying assumption fixes the
shape of $\varOmega$ in its referential configuration  also allows  for
considering the fixed boundary even for such time-evolving Eulerian description.
The other conditions in \eq{Euler-thermodynam-BC} model the free slip in the
tangential direction
in our multipolar variant
in  the  variationally legitimate simplest way, cf.\
\cite[Formula~(2.24)]{Roub24TVSE}. The last condition in
\eq{Euler-thermodynam-BC} prescribes the heat flux through the boundary.

The energetics behind the system is revealed when testing
\eq{Euler1-thermo-finite} by $\vv$ and \eq{Euler3-thermo-finite} by 1.
To execute the former test, let us use  \eqref{Euler0-thermo-finite}
tested by $\frac12|\vv|^2$
for $\pdt{}(\varrho|\vv|^2/2)=\varrho\vv{\cdot}\pdt{}\vv
-{\rm div}(\varrho\vv)|\vv|^2/2$ and realize that
$\int_\varOmega\varrho\DT\vv{\cdot}\vv\,\d\xx=
\frac{\d}{\d t}\int_\varOmega\frac12\varrho|\vv|^2\,\d\xx
+\int_\varGamma\varrho|\vv|^2\vv{\cdot}\nn\,\d S$. Thus, by using Green
formula twice relying also on \eq{Euler-thermodynam-BC}, we obtain
the {\it mechanical-energy dissipation balance}
\begin{align}
  &\!\!\!\!\!\frac{\d}{\d t}
  \int_\varOmega\!\!\!\!
  \linesunder{\frac\varrho2|\vv|^2}{kinetic}{energy}\!\!\!\!\d\xx
  +\int_\varOmega\!\!\!\!
\lineunder{\nu_1|\ee(\vv)|^2\!+\nu_2|\Nabla^2\vv|^p_{_{_{_{_{_{}}}}}}}
{dissipation rate $\xi$}\!\!\!\!\d\xx
=\int_\varOmega\!\!\!\!\linesunder{\varrho\GRAVITY{\cdot}\vv}{power of}{gravity}\!\!\!\!-\!\!\!\!\!\morelinesunder{\TT(\FF,\theta){:}\nabla\vv}{power of}{conservative}{stress}\!\!\!\!
\d\xx\,.
\label{Euler-engr-finite-}\end{align}
When added \eq{Euler3-thermo-finite} tested by 1,
we obtain the {\it total-energy balance}:
\begin{align}
  &\!\!\!\!\!\frac{\d}{\d t}
  \int_\varOmega\!\!\!\!
  \linesunder{\frac\varrho2|\vv|^2}{kinetic}{energy}\!\!\!\!+\!\!\!\!\!
  \linesunder{
\ENG(\FF,\theta)}{internal}{energy}\!\!\!\!\d\xx
=\int_\varOmega\!\!\!\!\!\!\!\!\linesunder{\varrho\GRAVITY{\cdot}\vv}{power of}{gravity field}\!\!\!\!\!\!\!\!\!
\d\xx
+\!\int_\varGamma\!\!\!\morelinesunder{h}{external}{heat}{flux}\!\!\!\d S\,,
\label{Euler-engr-finite}\end{align}
expressing the {\it 1st law of thermodynamics}: in isolated
systems, the total energy is conserved.

Furthermore, for the test of \eq{Euler3-thermo-finite} by the so-called coldness
$1/\theta$, we need to assume positivity of $\theta$ and then use the calculus
\begin{align}\label{calculus-to-entropy}
\frac1\theta\Big(\pdt\Eng+{\rm div}(\Eng\vv)\Big)
&=\frac1\theta\,\DT{{\overline{\ENG(\FF,\theta)}}}
+\frac{\ENG(\FF,\theta)}{\theta}{\rm div}\,\vv
\\&\nonumber=\frac{\ENG_\theta'(\FF,\theta)}{\theta}\DT\theta
+\frac{\ENG_\FF'(\FF,\theta)}{\theta}\Colon\DT\FF
+\frac{\ENG(\FF,\theta)}{\theta}{\rm div}\,\vv
\\&\nonumber=\DT{{\overline{\eta(\FF,\theta)}}}
+\Big(\frac{\ENG_\FF'(\FF,\theta)}{\theta}{-}\eta_\FF'(\FF,\theta)
\Big)\Colon\DT\FF+\frac{\ENG(\FF,\theta)}{\theta}{\rm div}\,\vv
\\&\nonumber=\DT{{\overline{\eta(\FF,\theta)}}}
+\eta(\FF,\theta){\rm div}\,\vv
+\frac{\psi_\FF'(\FF,\theta)}{\theta}\Colon\DT\FF
+\frac{\psi(\FF,\theta)}{\theta}{\rm div}\,\vv
\\&=\pdt{}\eta(\FF,\theta)+{\rm div}(\eta(\FF,\theta)\vv)
+\!\!\!\!\lineunder{\frac{\psi_\FF'(\FF,\theta)}{\theta}\Colon\DT\FF
+\frac{\psi(\FF,\theta)}{\theta}{\rm div}\,\vv}{$=\TT(\FF,\theta){:}\ee(\vv)/\theta$}\!\!\!\!\,.
\nonumber\end{align}
Here we also used $\eta_\theta'(\FF,\theta)=\ENG_\theta'(\FF,\theta)/\theta$ and
$\psi(\FF,\theta)/\theta-\ENG(\FF,\theta)/\theta=\psi_\theta'(\FF,\theta)=-\eta(\FF,\theta)$.
Therefore, by the mentioned test of \eq{Euler3-thermo-finite}, we obtain 
the {\it entropy balance}:
\begin{align}
  &\frac{\d}{\d t}\!\!\!\!
  \linesunder{\int_\varOmega\!\eta(\Fe,\theta)\,\d\xx}{total}{entropy}\!\!\!\!\!
  =\int_\varOmega\!\!\!\!\linesunder{\frac{
  \nu_1|\ee(\vv)|^2+\nu_2|\Nabla^2\vv|^p}\theta}{entropy production due to}{mechanical viscosity}\!\!\!\!+\!\!\!\!
  \morelinesunder{\kappa(\Fe,\theta)\frac{|\nabla\theta|^2}{\theta^2}}{entropy
  produ-}{ction due to}{heat transfer}\!\!\!\d\xx
+\int_\varGamma\!\!\!\!\!\!\!\!\!
\morelinesunder{\frac{h}\theta}{entropy flux}{through}{boundary}\!\!\!\!\!\!\!\!\!\d S\,,
  \label{Euler-entropy-finite}\end{align}
expressing the {\it 2nd law of thermodynamics}: in isolated systems, the total
entropy is not decaying.

Eliminating $\Eng=\ENG(\FF,\theta)=\psi(\FF,\theta)-\theta\psi_\theta'(\FF,\theta)$
from \eq{Euler3-thermo-finite} by using the calculus
\begin{align}\nonumber
\DT{\overline{\ENG(\FF,\theta)}}&=\psi_\FF'(\FF,\theta){:}\DT\FF+\psi_\theta'(\FF,\theta)\DT\theta-\DT\theta\psi_\theta'(\FF,\theta)
-\theta\psi_{\FF\theta}''(\FF,\theta)\DT\FF
-\theta\psi_{\theta\theta}''(\FF,\theta)\DT\theta
\\&\nonumber\!\!\!\stackrel{\scriptsize\eq{Euler2-thermo-finite}}{=}\!
\TT(\FF,\theta){:}\nabla\vv
-\theta\TT_\theta'(\FF,\theta){:}\nabla\vv-\Eng{\rm div}\,\vv\,,
\end{align}
we obtain the {\it heat-transfer equation} for
temperature
\begin{align}\label{Euler3-heat-eq}
\!\!c(\FF,\theta)\DT\theta={\rm div}\big(\kappa(\FF,\theta)\nabla\theta\big)
+\nu_1|\ee(\vv)|^2+\nu_2|\Nabla^2\vv|^p
+\theta\big(\psi_{\FF\theta}''(\FF,\theta)\FF^\top\!\!+\psi_\theta'(\FF,\theta)\bbI
\big){:}\nabla\vv\ \
\\\ \ \ \text{ with }\ \ c(\FF,\theta)=-\theta\psi_{\theta\theta}''(\FF,\theta)\,,
\nonumber\end{align}
where $c(\FF,\theta)$ is in the position of heat capacity. The heat equation
\eq{Euler3-heat-eq} involves explicitly $\psi_{\FF\theta}''$ and
$\psi_{\theta\theta}''$ and thus needs a  smoother  $\psi$ than
\eq{Euler3-thermo-finite}.

\begin{remark}[{\sl Selecting the mere stored temperature-independent energy out}]\label{rem-split}\upshape
Sometimes, it is useful to split the free energy as
$\psi(\FF,\theta)=\phi(\FF)+\varphi(\FF,\theta)$, where $\phi$ plays the role of
a stored energy. Such a splitting allows for a calibration $\varphi(\FF,0)=0$,
which actually determines this splitting uniquely since then
$\phi(\FF)=\psi(\FF,0)$ and will be used in Steps~3 and 4 in the proof of
Theorem~\ref{prop-Euler} below. This also implies the splitting of the Cauchy
stress as
\begin{align}\label{stress-split}
\TT(\FF,\theta)=\TT_0(\FF)+\TT_1(\FF,\theta)\ \ \ \text{ with }\ \ \TT_0(\FF):=\phi'(\FF)\FF^\top\!+\phi(\FF)\bbI
=\TT(\FF,0)\,,
\end{align}
which allows for usage of the calculus
\begin{align}\label{calcul-for-phi}
\TT_0(\FF){:}\nabla\vv
&=\phi'(\FF)\FF^\top{:}\nabla\vv+\phi(\FF){\rm div}\,\vv
=\phi'(\FF){:}((\nabla\vv)\FF)+\phi(\FF){\rm div}\,\vv
\\[-.3em]&\nonumber
\!\!\!\!\stackrel{\scriptsize\eq{Euler2-thermo-finite}}{=}\!\phi'(\FF){:}\DT\FF+\phi(\FF){\rm div}\,\vv
=\DT{\overline{\phi(\FF)}}+\phi(\FF){\rm div}\,\vv\,.
\end{align}
This further allows for the calculus
\begin{align}\nonumber
\int_\varOmega\!\TT_0(\FF){:}\nabla\vv\,\d\xx
&
=\frac{\d}{\d t}\int_\varOmega\!\phi(\FF)\,\d\xx
+\int_\varOmega\!\nabla\phi(\FF){\cdot}\vv+\phi(\FF){\rm div}\,\vv\,\d\xx=
\\&\nonumber
=\frac{\d}{\d t}\int_\varOmega\!\phi(\FF)\,\d\xx
+\!\int_\varOmega\!{\rm div}\big(\phi(\FF)\vv\big)\,\d\xx=
\frac{\d}{\d t}\int_\varOmega\!\phi(\FF)\,\d\xx+\!\int_\varGamma\!\phi(\FF)\vv{\cdot}\nn\,\d S\,,
\end{align}
where the last term vanishes due to the boundary condition $\vv{\cdot}\nn=0$ on
$\varGamma$.  In  this way, one can see a more detailed
mechanical-energy dissipation balance than \eq{Euler-engr-finite-}, namely 
\begin{align}
  &\!\!\!\!\!\frac{\d}{\d t}
  \int_\varOmega\!\!\!\!
  \linesunder{\frac\varrho2|\vv|^2}{kinetic}{energy}\!\!\!\!\!+\!\!\!\!
  \linesunder{\phi(\FF)}{stored}{energy}\!\!\!\d\xx
  +\int_\varOmega\!\!\!\!
\lineunder{\nu_1|\ee(\vv)|^2\!
+\nu_2|\Nabla^2\vv|^p_{_{_{_{_{_{}}}}}}}
{dissipation rate}\!\!\!\!\d\xx
=\int_\varOmega\!\!\!\!\linesunder{\varrho\GRAVITY{\cdot}\vv}{power of}{gravity}\!\!\!\!
-\!\!\!\!\!\morelinesunder{\TT_1(\FF,\theta){:}\nabla\vv}{power of}{adiabatic}{effects}\!\!\!\!
\d\xx\,.
\label{Euler-engr-finite--}\end{align}
Also the internal energy $\ENG$ used in the total energy balance
\eq{Euler-engr-finite} then sees directly the stored energy, since obviously
$\ENG(\FF,\theta)=\phi(\FF)+\varphi(\FF,\theta)-\theta\varphi_\theta'(\FF,\theta)$
and the mentioned calibration $\varphi(\FF,0)=0$ yields that simply
$\phi(\FF)=\ENG(\FF,0)=\psi(\FF,0)$. Noteworthy, using \eq{calcul-for-phi} for
the internal-energy equation \eq{Euler3-thermo-finite} eliminates the mere
stored energy and yields some {\it heat-type equation} in an enthalpy-like
formulation
\begin{align}\label{Euler3-enthalpy-finite}
&
\DT\Ent={\rm div}\big(\kappa(\FF,\theta)\nabla\theta\big)+\nu_1|\ee(\vv)|^2
+\nu_2|\Nabla^2\vv|^p
-({\rm div}\,\vv)\Ent+\TT_1(\FF,\theta){:}\nabla\vv
\\&\hspace{3em}\text{ with }\ \ \Ent=\ENT(\FF,\theta)\,,
\ \text{ where }\ \ENT(\FF,\theta):=
\psi(\FF,\theta)-\theta\psi_\theta'(\FF,\theta)-\psi(\FF,0)\,.
\nonumber
\end{align}
Here $\TT_1(\FF,\theta)=\TT(\FF,\theta)-\TT(\FF,0)$ is from \eq{stress-split}.
Note that $\ENT(\FF,0)=0$ and $\TT_1(\FF,0)=\bm0$.
\end{remark}

\begin{remark}[{\sl Thermoelasticity in an ``engineering'' formulation}]
\label{rem-engineering}\upshape
The system \eq{Euler-thermo-finite} is well prepared for analysis but
may not be found ``engineering friendly''. There are two aspects:
the free energy is more often considered referential rather than actual
(let us distinguish it by $\uppsi$) and the heat equation is formulated
in temperature rather than internal energy as in \eq{Euler3-thermo-finite}.
The free energy by reference volume is more standard in
continuum physics \cite{GuFrAn10MTC,Mart19PCM} than the free energy per
actual evolving volume, as it corresponds more directly to
experimentally available data. The relation in between
$\uppsi$ and $\psi$ is $\psi(\FF,\theta)=\uppsi(\FF,\theta)/\!\det\FF$ and
then the conservative part of the (actual) Cauchy stress, the 
actual entropy, and the actual heat capacity are respectively
\begin{align}
\TT(\FF,\theta)=\frac{\uppsi_\Fe'(\Fe,\theta)\Fe^\top\!\!}{\det\Fe}\,,\ \ \ \ \ 
\eta(\FF,\theta)=-\,\frac{\!\uppsi_\theta'(\Fe,\theta)\!}{\det\Fe}\,,
\ \ \text{ and }\ \ 
c(\FF,\theta)=-\theta\,\frac{\!\uppsi_{\theta\theta}''(\FF,\theta)\!}{\det\FF}\,.
\label{stress-entropy-eng}\end{align}
The resulting system then consists from (\ref{Euler-thermo-finite}a-c)
but with $\TT$ from \eq{stress-entropy-eng}. Moreover, \eq{Euler3-thermo-finite} 
in the form of the heat equation like \eq{Euler3-heat-eq} reads as
\begin{align}
&c(\FF,\theta)\DT\theta={\rm div}\big(\kappa(\FF,\theta)\nabla\theta\big)
+\nu_1|\ee(\vv)|^2\!+\nu_2|\nabla^2\vv|^p\!
+\theta\frac{\uppsi_{\FF\theta}''(\FF,\theta)\FF^\top\!\!\!}{\det\FF}{:}\nabla\vv\,.
\label{heat-eq+}\end{align}
Actually, the split  $\TT=\TT_0+\TT_1$  from Remark~\ref{rem-split}
allows for a specification of the last term 
 showing that  the mere stored energy  does not contribute 
to the adiabatic power  since $[\TT_0]_\theta'\equiv0$ so that this  last
term in \eq{heat-eq+} can
 actually  be written as $\theta[\TT_1]_\theta'(\FF,\theta){:}\nabla\vv$.
 Moreover,  in literature the factor $1/\det\FF$ is sometimes replaced
by $\varrho$ if the referential free energy $\uppsi$ is  in  the physical
units  J/kg  instead of  Pa=J/m$^3$.
\end{remark}

\section{Analytical results: existence of weak solutions}\label{sec-anal}

To highlight the nontriviality of the coupled system \eq{Euler-thermo-finite}, 
it is perhaps useful to give an overview of various strategies in the proof
that can or should be considered. While the momentum equation is to be
ultimately tested by the velocity $\vv$, the heat equation can be subjected
to various tests by:\\
\Item{\ \,(i)}{1 to obtain energy balance,}
\Item{\ (ii)}{$1/\theta$ to obtain entropy balance (not directly
exploited in this paper, see Remark~\ref{rem-entrop}),}
\Item{\,(iii)}{$\theta^-$ to see non-negativity of temperature (not used in
this paper),}
\Item{\,(iv)}{$\theta$ to see an estimate for $\nabla\theta$  (which
will need a regularization of the heat
sources), and}
\Item{\ (v)}{$1-1/(1{+}\theta)^\ZETA$ for $\ZETA>0$ to estimate $\nabla\theta$
for the physically relevant $L^1$-heat sources.}
In (iii), we have used the notation (decomposition)
\begin{align}
\theta=\theta^+-\theta^-\ \ \ \text{ with }\ \theta^+:=\max(0,\theta)\,.
\end{align}
The construction of an approximate solution by time-discretization is very
problematic because the free energy $\psi(\cdot,\theta)$ is highly nonconvex.
Also the Galerkin discretization is nontrivial even if applied
only to the momentum and the heat equations while keeping the
transport equations continuous, as used below. The nontriviality of
this scenario consists in the mentioned tests of the heat equation:
only the tests (i) and (iv) are legitimate in conformal Galerkin approximations
and, for (i), the internal energy is to be bounded from below
while, for (iv), the heat capacity should be granted non-negative for negative
temperatures.  These requirements are  slightly
contradictory:  the boundedness of the internal energy $\ENG$ from below
essentially means that $\ENG$ 
stays constant and the heat capacity $c$ simply vanishes for negative
temperatures,  which however is not suitable for the test (iv) directly.
For this reason, we devise a suitable extension of data for negative
temperature values,
cf.\ the
illustration in Figure~\ref{various-nonlinearities},
 which will give in the limit  the non-negativity of temperature as a side
effect without using the test (iii).

We will use the standard notation concerning the Lebesgue and the Sobolev
spaces, namely $L^p(\varOmega;\R^n)$ for Lebesgue measurable functions
$\varOmega\to\R^n$ whose Euclidean norm is integrable with $p$-power, and
$W^{k,p}(\varOmega;\R^n)$ denotes  the space of  functions from
$L^p(\varOmega;\R^n)$ whose all derivatives up to the order $k$ have their
Euclidean norm integrable with $p$-power. We also write briefly $H^k=W^{k,2}$.
The notation $p'$  will denote the conjugate  exponent $p/(p{-}1)$ while  
$p^*$ will denote the exponent from the embedding
$W^{1,p}(\varOmega)\subset L^{p^*}(\varOmega)$, i.e.\ $p^*=dp/(d{-}p)$
for $p<d$ while $p^*\ge1$ is arbitrary for $p=d$ or $p^*=+\infty$ for $p>d$.
Furthermore, for a Banach space
$X$ and for $I=[0,T]$, we will use the notation $L^p(I;X)$ for the Bochner
space of Bochner measurable functions $I\to X$ whose norm is in $L^p(I)$
while $W^{1,p}(I;X)$ denotes for functions $I\to X$ whose distributional
derivative is in $L^p(I;X)$. Also, $C(\cdot)$ and $C^1(\cdot)$
will denote the spaces of continuous and continuously differentiable functions,
respectively, and $\barOmega$ will denote the closure of $\varOmega$. Eventually,
the spaces of weakly continuous functions $I\to X$ is denoted by $C_{\rm w}(I;X)$.

Moreover, as usual, we will use $C$ for a generic constant which may vary
from estimate to estimate  whose value depends on the data qualified in
\eq{Euler-ass} below. Occasionally, when such constants depend also on some
discretization/regularization parameters during the particular steps in the
proof below, we will indicate this dependence by a lower index.

We will consider an initial-value problem, prescribing the initial conditions
\begin{align}\label{Euler-thermodynam-IC}
\varrho|_{t=0}^{}=\varrho_0:=\frac{\rhoR}{\det\FF_0}\,,
\ \ \ \ \ \vv|_{t=0}^{}=\vv_0\,,\ \ \ \ \ 
\FF|_{t=0}^{}=\FF_0\,,\ \ \text{ and }\ \ \theta|_{t=0}^{}=\theta_0\,.
\end{align}

To devise a weak formulation of the initial-boundary-value problem
\eq{Euler-thermodynam-BC} and \eq{Euler-thermodynam-IC} for
the system \eq{Euler-thermo-finite}, we use the by-part integration in time and
the Green twice formula for \eq{Euler1-thermo-finite} multiplied by
a smooth test function $\widetilde\vv$ and once for \eq{Euler3-thermo-finite}
multiplied by a smooth test function $\widetilde\theta$, while
(\ref{Euler-thermodynam-BC}a,c) are considered just a.e.\ on $I{\times}\varOmega$:

\begin{definition}[Weak solutions to \eq{Euler-thermo-finite}]\label{def}
For $p,q\in[1,\infty)$, a quadruple $(\varrho,\vv,\Fe,\theta)$ with
$\varrho\in L^\infty(I{\times}\varOmega)\cap W^{1,1}(I{\times}\varOmega)$,
$\vv\in L^2(I;H^1(\varOmega;\R^d))\,\cap\,L^p(I;W^{2,p}(\varOmega;\R^d))$,
$\FF\in L^\infty(I{\times}\varOmega;\R^{d\times d})
\cap W^{1,1}(I{\times}\varOmega;\R^{d\times d})$,
and $\theta\in L^1(I;W^{1,1}(\varOmega))$ will be called a weak
solution to the system \eq{Euler-thermo-finite} with the boundary conditions
\eq{Euler-thermodynam-BC} and the initial condition \eq{Euler-thermodynam-IC}
if the integral identities 
\begin{subequations}\label{Euler-weak}\begin{align}
&
\int_0^T\!\!\!\!\int_\varOmega\bigg(\Big(
\psi_\Fe'(\Fe,\theta)\Fe^\top\!\!+
\nu_1\ee(\vv)-\varrho\vv{\otimes}\vv\Big){:}\ee(\widetilde\vv)
+\psi(\Fe,\theta){\rm div}\,\widetilde\vv
\label{Euler1-weak}\\[-.4em]&\hspace{1.5em}
+\nu_2|\nabla^2\vv|^{p-2}\nabla^2\vv\Vdots\nabla^2\widetilde\vv
-\varrho\vv{\cdot}\pdt{\widetilde\vv}\bigg)\,\d\xx\d t
=\!\int_0^T\!\!\!\!\int_\varOmega\varrho\GRAVITY{\cdot}\widetilde\vv\,\d\xx\d t
+\!\int_\varOmega\!\varrho_0\vv_0{\cdot}\widetilde\vv(0)\,\d\xx
\ \text{ and}
\nonumber\\\label{Euler3-weak}
&\int_0^T\!\!\!\int_\varOmega\!\bigg(\!\ENG(\FF,\theta)\pdt{\wt\theta}
+\Big(\ENG(\FF,\theta)\vv{+}\kappa(\FF,\theta)\nabla\theta\Big){\cdot}\nabla\wt\theta
+\Big(\nu_1|\ee(\vv)|^2+\nu_2|\nabla^2\vv|^p\!
\\[-.4em]&\hspace{6em}
+\TT(\FF,\theta){:}\nabla\vv\!\Big)\wt\theta\bigg)\d\xx\d t
+\!\int_0^T\!\!\!\int_\varGamma h\wt\theta\,\d S\d t
+\!\int_\varOmega\!\ENG(\FF_0,\theta_0)\wt\theta(0)\,\d\xx=0
\nonumber
\end{align}
\end{subequations}
with $\TT(\cdot,\cdot)$  from \eq{Euler1-thermo-finite} and $\ENG(\cdot,\cdot)$
from \eq{Euler3-thermo-finite} hold
for any $\widetilde\vv$ smooth with $\widetilde\vv{\cdot}\nn={\bm0}$ and
$\widetilde\vv(T)=0$ and for any $\wt\theta$ smooth with
$\wt\theta(T)=0$, and eventually if also \eq{Euler0-thermo-finite} and
\eq{Euler2-thermo-finite} hold a.e.\ on $I{\times}\varOmega$ together with the
initial conditions for $\varrho$ and $\FF$ in \eq{Euler-thermodynam-IC}.
\end{definition}

Before stating the main analytical result, let us summarize the data
qualification. We use the notation $\R^+:=[0,+\infty)$. For some $\delta>0$, and
some $1<q<\infty$ and $d<p<\infty$, recalling the notation 
\begin{align}
\ENG(\FF,\theta)=\psi(\FF,\theta)-\theta\psi_\theta'(\FF,\theta)
\ \ \ \text{ and }\ \ \
\TT(\FF,\theta)=\psi_\FF'(\FF,\theta)\FF^\top\!+\psi(\FF,\theta)\bbI\,,
\end{align}
and denoting by ${\rm GL}^+(d)=\{F\in\R^{d\times d};\ \det F>0\}$ the
orientation-preserving general linear group, we assume:
\begin{subequations}\label{Euler-ass}\begin{align}
&\varOmega\ \text{ a smooth bounded domain of $\R^d$, }\ d=2,3,
\\&\label{Euler-ass-psi}
\psi\in C^2({\rm GL}^+(d)\times\R)\,,\ \ \kappa
\in C({\rm GL}^+(d){\times}\R)\,,
\\[-.1em]&\nonumber
\exists\,C\in\R\ \ \forall(\FF,\theta)\in{\rm GL}^+(d){\times}\R:
\\[-.1em]&\hspace{2em}
\psi(\FF,\theta^-)=\ENG(\FF,0)\ \ \text{ and }\ \ \kappa(\FF,\theta^-)=\kappa(\FF,0)\,,
\label{Euler-ass-psi-increas}
\\[-.1em]&\hspace{2em}
\label{ass-stress-control}
|\TT(\FF,\theta)|\le C\big(1\,{+}\,\ENG(\FF,\theta)\big)\,,
\\[-.6em]&\nonumber
\forall K\!\subset\!{\rm GL}^{\!+}(d)\text{ compact }
\ \,\exists\, 0{<}c_K^{}{\le} C_K^{}{<}+\infty,\ 
0\le\ALPHA<
\begin{cases}\ 1&\text{for }d=2\\[-.3em]1/2\!\!\!\!&\text{for }d=3\end{cases}
\ \forall(\FF,\theta){\in}K{\times}\R{:}
\\[-.2em]&\label{Euler-ass-psi-1}
\hspace{2em}
\ENG(\FF,\theta)\ge  
c_K^{}(\theta^+)^{1+\ALPHA}\ \,\text{ and }\ \ \ENG_\theta'(\FF,\cdot)>0\text{ 
on $(0,+\infty)$\,,}
\\[-.1em]&\hspace{2em}\label{Euler-ass-psi-1+++}
|\ENG_\FF'(\FF,\theta)|
+\theta^+\ENG_\theta'(\FF,\theta)
\le C_K^{}\big(1\,{+}\,(\theta^+)^{1+\ALPHA}\big)\ \text{ and }\ 
|\ENG_{\FF\theta}''(\FF,\theta)|\le C_K^{}\big(1\,{+}\,(\theta^+)^{\ALPHA}\big),
\\&\hspace{2em}
c_K^{}\le\kappa(\FF,\theta)\le C_K^{}\,,\ \ \nu_1>0,\ \nu_2>0\,,
\label{Euler-ass-kappa}
\\&\vv_0\in L^2(\varOmega;\R^d)\,,\ \ \ \
\Fe_0\in W^{1,\infty}(\varOmega;\R^{d\times d})\ \ \text{ with }\ \
{\rm min}_{\barOmega}^{}\det\Fe_0>0\,,
\label{Euler-ass-Fe0}
\\[-.1em]&
{\bm g}\in L^1(I;L^\infty(\varOmega;\R^d))\,,\,\ \ 
\rhoR\in W^{1,\infty}(\varOmega)
\ \ \text{ with }\ \ \ \ {\rm min}_{\barOmega}^{}\rhoR>0\,,
\label{Euler-ass-rhoR}
\\[-.1em]&h\,{\in}\, L^1(I{\times}\varGamma)\,,\ \ h\ge0\,,\ \ \theta_0\in L^1(\varOmega),\ \ \ \theta_0\ge0\ \text{ a.e.\ on }\ \varOmega\,,\ \ \ 
\ENG(\FF_0,\theta_0)\in L^1(\varOmega)\,.
\label{Euler-ass-theta0}
\end{align}\end{subequations}
Let us note that we have formally defined the data also for negative
temperatures, which will be useful in the proof below.
Noting $\ENG_\theta'(\Fe,\theta)=-\theta\psi_{\theta\theta}''(\Fe,\theta)$,
the latter condition in \eq{Euler-ass-psi-1} implies $\psi(\FF,\cdot)$
concave on $\R^+$. The relevance of the possibility of controlling the stress
via the energy, which is what \eq{ass-stress-control} does,
was pointed out by J.M.\,Ball \cite{Ball84MELE,Ball02SOPE}
for the Kirchhoff rather than the Cauchy stress, cf.\ Remark~\ref{rem-control}
below. Actually, \eq{Euler-ass-psi} allows also for  
$\sup_{\FF\in K}|\psi(\FF,0)|\le c_K$ for some $C_K$ since
$\ENG(\cdot,0)=\psi(\cdot,0)\in C({\rm GL}^+(d))$, and then
\eq{Euler-ass-psi-1+++} gives, for $K$ and $\ALPHA$ as in
(\ref{Euler-ass}e--g), a bound for the internal energy $\ENG$ as
\begin{align}\label{growth-int-eng}
\forall (\FF,\theta)\in K{\times}\R:\ \ \ \ENG(\FF,\theta)=\psi(\FF,0)+
\!\int_0^\theta\!\!\ENG_\theta'(\FF,\wt\theta)\,\d\wt\theta\le
C_K\big(1+(\theta^+)^{1+\ALPHA}\big)
\end{align}
so that the energy-controlled-stress condition \eq{ass-stress-control}
gives
\begin{align}\label{growth-T}
\forall (\FF,\theta)\in K{\times}\R:\ \ \ 
|\TT(\FF,\theta)|\le C_K\big(1+(\theta^+)^{1+\ALPHA}\big)\,.
\end{align}
Moreover, the expected {\it symmetry} of such Cauchy stress $\TT$
is granted by {\it frame indifference} of $\psi(\cdot,\theta)$. This means
that 
\begin{align}
\forall (\FF,\theta)\in{\rm GL}^+(d){\times}\R,\ \ Q\in{\rm SO}(d):\ \ \ \ 
\psi(\FF,\theta)=\psi(Q\FF,\theta)\,,
\label{frame-indifference}\end{align}
where $Q\in{\rm SO}(d)=\{Q\in\R^{d\times d};\ Q^\top Q=QQ^\top=\bbI\}$ is the
special orthogonal group.

\begin{theorem}[Existence and regularity of weak solutions]\label{prop-Euler}
Let $p>d$ and the assumptions \eq{Euler-ass} and \eq{frame-indifference} hold.
Then:\\
\Item{(i)}{there exists at least one weak solution $(\varrho,\vv,\FF,\theta)$ according
Definition~\ref{def} such that, in addition,
$\varrho\in C_{\rm w}(I;W^{1,r}(\varOmega))$ and
$\Fe\in C_{\rm w}(I;W^{1,r}(\varOmega;\R^{d\times d}))$ for any $1\le r<\infty$,
and further $\theta\in C_{\rm w}(I;L^{1+\ALPHA}(\varOmega))
\,\cap\,L^\EXP(I;W^{1,\EXP}(\varOmega))$ with $1\le\EXP<(d{+}2{+}(2{-}d)\ALPHA)/(d{+}1{+}\ALPHA)$.
Moreover,
$\pdt{}(\varrho\vv)\in L^2(I;H^1(\varOmega;\R^d)^*)+L^{p'}(I;W^{2,p}(\varOmega;\R^d)^*)$
and $\,\inf_{I{\times}\varOmega}^{}\varrho>0$, $\,\inf_{I{\times}\varOmega}^{}\det\FF>0$,
and $\theta\ge0$ a.e.\ on $I{\times}\varOmega$.}
\Item{(ii)}{This  solution complies with energetics in the sense that
the energy-dissipation balances  \eq{Euler-engr-finite-} and
\eq{Euler-engr-finite--} as well as the total energy balance
\eq{Euler-engr-finite} integrated over time interval $[0,t]$ with the
initial conditions \eq{Euler-thermodynam-IC} hold.}
\end{theorem}

\begin{proof}
For clarity, we will divide the proof into  seven  steps.

\medskip\noindent{\it Step 1: Regularization of the system and a Galerkin
semi-discretization}.
For $\EPS>0$ and $k\in\N$, we consider a regularization of the momentum
equation \eq{Euler1-thermo-finite} and of the internal-energy equation
\eq{Euler3-thermo-finite}. Altogether, the system \eq{Euler-thermo-finite}
is regularized as
\begin{subequations}\label{Euler-thermo-finite+}\begin{align}
&\DT\varrho=-({\rm div}\,\vv)\varrho\,,
\label{Euler0-thermo-finite+}
\\\label{Euler1-thermo-finite+}&\varrho\DT\vv={\rm div}\big(\TT(\FF,\theta){+}\Sv\big)
 +\varrho\GRAVITY-\frac1k\vv\ \ \
\text{ where }\
\TT(\FF,\theta)=\psi_\FF'(\FF,\theta)\FF^\top\!\!+\psi(\FF,\theta)\bbI
\\[-.4em]&\hspace*{16em}\text{ and }\ \ 
\Sv=\nu_1\ee(\vv)-{\rm div}\big(\nu_2|\nabla^2\vv|^{p-2}\nabla^2\vv\big)\,,
\nonumber
\\\label{Euler2-thermo-finite+}
&\DT\FF=(\nabla\vv)\FF\,,
\\[-.3em]
&\label{Euler3-thermo-finite+}
\DT\Eng={\rm div}\big(\kappa(\FF,\theta)\nabla\theta\big)-({\rm div}\,\vv)\Eng
+\frac{\nu_1|\ee(\vv)|^2+\nu_2|\Nabla^2\vv|^p}{1{+}\EPS|\ee(\vv)|^2{+}\EPS|\Nabla^2\vv|^p}
+\TT(\FF,\theta){:}\nabla\vv-\frac1k\pdt\theta
\\[-.1em]&\hspace{18.5em}\ \ \ \text{ with }\ \ \Eng=\psi(\FF,\theta)
-\theta\psi_\theta'(\FF,\theta)\,.
\nonumber\end{align}\end{subequations}
The equations (\ref{Euler-thermo-finite+}b,d) are completed by the 
boundary conditions \eq{Euler-thermodynam-BC} but with
the last initial condition regularized, specifically 
\begin{align}\label{Euler-thermodynam-BC-2-reg}
\kappa(\FF,\theta)\nabla\theta{\cdot}\nn=h_\EPS\ \ \ \text{ with }\ \
h_\EPS:=\frac{h}{1{+}\EPS h}
\end{align}
and similarly the initial conditions \eq{Euler-thermodynam-IC} are considered
with a regularized temperature, namely
\begin{align}\label{Euler-thermodynam-IC-reg}
\varrho|_{t=0}^{}=\varrho_0\,,\ \ \ \ \ \vv|_{t=0}^{}=\vv_0\,,\ \ \ \ \ 
\FF|_{t=0}^{}=\FF_0\,,\ \ \text{ and }\ \ \theta|_{t=0}^{}=\theta_{0\EPS}:=
\frac{\theta_0}{1{+}\EPS\theta_0}\,.
\end{align}
The purpose of the regularization $ k^{-1}\vv$ in
\eq{Euler1-thermo-finite+} is to make values of $\vv$ controlled
 (even not uniformly with respect to $k$) already  before we
 will  establish a uniform estimate about
the sparsity $1/\varrho$  and then $\vv$ in (\ref{est+}b,c);
cf.\ the proof in \cite[Lemma~5.1]{Roub24TVSE}, while the $\EPS$-regularization
in \eq{Euler3-thermo-finite+} and in \eq{Euler-thermodynam-BC-2-reg} is to allow
its testing by $\theta$.  Specifically,  the regularizing term
$ k^{-1}\pdt{}\theta$ in
\eq{Euler3-thermo-finite+} will ensure the existence of  a global-in-time 
Galerkin solution by the usual successive-prolongation argument in Step~3 below
 exploiting the $L^\infty(I)$-estimates derived further in Steps~2 and 3.

We use a spatial semi-discretization, keeping the transport equations
\eq{Euler0-thermo-finite+} and \eq{Euler2-thermo-finite+} continuous
(i.e.\ non-discretised). More specifically, we make a conformal Galerkin
approximation of \eq{Euler1-thermo-finite+} by using a collection of 
nested finite-dimensional subspaces $\{V_k\}_{k=0,1,...}$ whose union is dense in
$W^{2,p}(\varOmega;\R^d)$ and a conformal Galerkin approximation of
\eq{Euler3-thermo-finite+} by using a collection of nested
finite-dimensional subspaces $\{Z_k\}_{k=0,1,...}$ whose union is dense
in $H^1(\varOmega)$. Note that we use the index $k$ occuring also
in \eq{Euler1-thermo-finite+} and in \eq{Euler3-thermo-finite+}.

Let us denote the solution of thus approximated regularized system
\eq{Euler-thermo-finite+} with the initial/boundary conditions
\eq{Euler-thermodynam-IC-reg} and \eq{Euler-thermodynam-BC} regularized as
\eq{Euler-thermodynam-BC-2-reg} by $(\varrho_{\EPS k},\vv_{\EPS k},\FFepsk,\theta_{\EPS k}):
I\to W^{1,r}(\varOmega)\times V_k\times W^{1,r}(\varOmega;\R^{d\times d})\times Z_k$.
Without loss of generality, we can assume $\vv_0\in V_0$ and $\theta_{0\EPS}\in Z_0$.
Specifically, $\varrho_{\EPS k}\in C_{\rm w}(I;W^{1,r}(\varOmega))\,\cap\,
W^{1,1}(I;L^1(\varOmega))$ and $\FFepsk\in C_{\rm w}(I;W^{1,r}(\varOmega;\R^{d\times d}))\,\cap\,
W^{1,1}(I;L^1(\varOmega;\R^{d\times d}))$ should satisfy
\begin{align}\label{transport-equations}
\pdt{\varrho_{\EPS k}}=-{\rm div}(\varrho_{\EPS k}\vv_{\EPS k})\ \ \ \text{ and }
\ \ \ \pdt{\FFepsk}=(\nabla\vv_{\EPS k})\FFepsk
-(\vv_{\EPS k}{\cdot}\nabla)\FFepsk\ \ \text{ a.e.\ on $I{\times}\varOmega$}
\end{align}
together with the following integral identities
\begin{subequations}\label{Euler-weak-Galerkin}\begin{align}
&\label{Euler1-weak-Galerkin}
\int_0^T\!\!\!\int_\varOmega
\bigg(\Big(\psi_\FF'(\FFepsk,\theta_{\EPS k})\FFepsk^\top\!
+\nu_1\ee(\vv_{\EPS k})
-\varrho_{\EPS k}\vv_{\EPS k}{\otimes}\vv_{\EPS k}\Big){:}\ee(\widetilde\vv)
\\[-.4em]&\hspace{1em}\nonumber
+\psi(\FFepsk,\theta_{\EPS k}){\rm div}\widetilde\vv
-\varrho_{\EPS k}\vv_{\EPS k}{\cdot}\pdt{\widetilde\vv}
+\nu_2|\nabla^2\vv_{\EPS k}|^{p-2}\nabla^2\vv_{\EPS k}\Vdots
\Nabla^2\widetilde\vv
+\frac{\!\vv_{\EPS k}{\cdot}\widetilde\vv}{k}\bigg)\,\d\xx\d t
\\[-.4em]&\hspace{16em}
=\!\int_0^T\!\!\!\int_\varOmega\varrho_{\EPS k}\GRAVITY{\cdot}\widetilde\vv\,\d\xx\d t
+\!\int_\varOmega\!\varrho_0\vv_0{\cdot}\widetilde\vv(0)\,\d\xx
\nonumber
\intertext{for any $\widetilde\vv\in L^\infty(I;V_k)$ with
$\widetilde\vv{\cdot}\nn=0$ on $I{\times}\varGamma$
and $\widetilde\vv(T)={\bm0}$, and}
\label{Euler3-weak-Galerkin}
&\!\int_0^T\!\!\!\int_\varOmega\bigg(\Big(\Eng_{\EPS k}{+}\frac{\theta_{\EPS k}}k\Big)\pdt{\widetilde\theta}
+\big(\Eng_{\EPS k}\vv_{\EPS k}
{-}\kappa(\FFepsk,\theta_{\EPS k})\nabla\theta_{\EPS k}\big)
{\cdot}\nabla\widetilde\theta
\\[-.3em]&\hspace{2em}\nonumber
+\Big(\frac{\nu_1|\ee(\vv_{\EPS k})|^2+\nu_2|\nabla^2\vv_{\EPS k}|^p}
{1{+}\EPS|\ee(\vv_{\EPS k})|^2{+}\EPS|\nabla^2\vv_{\EPS k}|^p}
+\TT(\FFepsk,\theta_{\EPS k}){:}\nabla\vv_{\EPS k}\Big)\widetilde\theta\,\bigg)
\d\xx\d t
\\[-.1em]&\hspace{4em}
+\!\int_\varOmega\!\ENG(\FF_0,\theta_{0\EPS})\widetilde\theta(0)\,\d\xx
+\!\int_0^T\!\!\!\int_\varGamma\!h_\EPS\widetilde\theta\,\d S\d t=0
\ \ \text{ with }\  \Eng_{\EPS k}=\ENG(\FFepsk,\theta_{\EPS k})
\nonumber
\end{align}
\end{subequations}
holding for any $\widetilde\theta\in C^1(I;Z_k)$ with
$\widetilde\theta(T)=0$, considered for $k=1,2,...$

 As announced in the previous,  existence of this solution is based on
the standard theory of systems of  abstract 
ordinary differential equations  valued in Banach spaces, here in
$W^{1,r}(\varOmega)\times V_k\times W^{1,r}(\varOmega;\R^{d\times d})\times Z_k$ 
first locally in time  by using standard arguments of the Cauchy–Peano theorem
 and then by successive prolongation on the whole time interval based on the
$L^\infty$-estimates below in Steps~2 and 3,  specifically \eq{est+} and
\eq{est-theta+}. For the  
the abstract $W^{1,r}(\varOmega)$-valued differential equations
\eq{transport-equations}  cf.\ also  \cite[Sect.5]{Roub24TVSE} or, in a bit
less general form, also \cite{Roub22QHLS,RouSte22VESS}.

\medskip\noindent{\it Step 2: first a-priori estimates}.
Actually, the fixing initial conditions $\FF_0\in W^{1,\infty}(\varOmega;\R^{d\times d})$
and $\varrho_0\in W^{1,\infty}(\varOmega)$ and $\vv_{\EPS k}$ enough regular, the
essential point is that the transport equations \eq{transport-equations} have
unique solutions. According \cite[Sect.5]{Roub24TVSE}, it thus defines the
weakly-continuous nonlinear operators
$\mathfrak{F}:I\times L^p(I;W^{2,p}(\varOmega;\R^d))\to
W^{1,r}(\varOmega;\R^{d\times d})$ and $\mathfrak{R}:I\times L^p(I;W^{2,p}(\varOmega;\R^d))
\to W^{1,r}(\varOmega)$ by 
\begin{align}\label{transport-equations+}
\FFepsk(t)=\mathfrak{F}\big(t,\vv_{\EPS k}\big)\ \ \text{ and }\ \
\varrho_{\EPS k}(t)=\mathfrak{R}\big(t,\vv_{\EPS k}\big)\,.
\end{align}
Moreover,  these nonlinear operators are bounded in the sense that, for any
$1\le r<+\infty$, $\mathfrak{F}(\cdot,\vv_{\EPS k})$ is bounded in 
$L^\infty(I;W^{1,r}(\varOmega;\R^{d\times d}))\,\cap\,
W^{1,1}(I;L^r(\varOmega;\R^{d\times d}))$ and $\mathfrak{R}(\cdot,\vv_{\EPS k})$ is
bounded in $L^\infty(I;W^{1,r}(\varOmega))\,\cap\,W^{1,1}(I;L^r(\varOmega))$
when $\nabla\vv_{\EPS k}$ ranges over a bounded set in
$L^p(I;W^{1,p}(\varOmega;\R^{d\times d}))$ provided also $\vv_{\EPS k}$  ranges
over a bounded set in  $L^1(I;W^{2,p}(\varOmega;\R^d))$ with
$\vv_{\EPS k}{\cdot}\nn=0$. Besides, the {\it mass conservation} holds:
\begin{align}\label{mass-conservation}
\varrho_{\EPS k}\ge0\ \text{ on }\ I{\times}\varOmega\ \ \ \text{ and }\ \ \ 
\int_\varOmega\!\varrho_{\EPS k}(t)\,\d\xx
=\!\int_\varOmega\!\varrho_0\,\d\xx\ \text{ for all }\ t\in I\,.
\end{align}

Testing \eq{Euler1-thermo-finite+} discretized as \eq{Euler1-weak-Galerkin}
by $\wt\vv=\vv_{\EPS k}$ gives the discrete mechanical-energy balance
as an analog of \eq{Euler-engr-finite-}:
\begin{align}\label{Euler-mech-engr-finite+}
\frac{\d}{\d t}\int_\varOmega\!\frac{\varrho_{\EPS k}}2|\vv_{\EPS k}|^2\,\d\xx+
  \!\int_\varOmega\!\nu_1|\ee(\vv_{\EPS k})|^2\!
+\nu_2|\nabla^2\vv_{\EPS k}|^p\!+\frac1k|\vv_{\EPS k}|^2\,\d\xx\ 
\\[-.4em]
=\int_\varOmega\varrho_{\EPS k}\GRAVITY{\cdot}\vv_{\EPS k}
-\TT(\FFepsk,\theta_{\EPS k}){:}\nabla\vv_{\EPS k}\,\d\xx\,.
\nonumber\end{align}
Testing \eq{Euler3-thermo-finite+} discretized as \eq{Euler3-weak-Galerkin} by
1 and summing it with \eq{Euler-mech-engr-finite+}, we obtain the variant of the
total-energy balance \eq{Euler-engr-finite} as an inequality:
\begin{align}
  &\frac{\d}{\d t}
  \int_\varOmega\!\frac{\varrho_{\EPS k}}2|\vv_{\EPS k}|^2+\ENG(\FFepsk,\theta_{\EPS k})
  \,\d\xx+\!\int_\varOmega\frac1k|\vv_{\EPS k}|^2\,\d\xx
\le\int_\varOmega\varrho_{\EPS k}\GRAVITY{\cdot}\vv_{\EPS k}\,\d\xx
+\!\int_\varGamma h_\EPS\,\d S\,.
\label{Euler-engr-finite+}\end{align}
The inequality (instead of an equality) in
\eq{Euler-engr-finite+} arises from the regularization of the
dissipative heat in \eq{Euler3-thermo-finite+}.

Moreover, using the calculus for the Galerkin solution \eq{calcul-for-phi}
 and introducing the notation for  the convective time derivative
related to the velocity field $\vv_{\EPS k}$ as
\begin{align}\label{def-of-DTek}
\DTk{\Ent\,}:=\pdt\Ent+\vv_{\EPS k}{\cdot}\Nabla\Ent\,,
\end{align}
the internal-energy equation \eq{Euler3-thermo-finite+} can be modified
to a heat equation  of the type  \eq{Euler3-enthalpy-finite}, i.e.\ here:
\begin{align}\label{Euler3-enthalpy-finite-reg}
\DTk\Ent_{\EPS k}&
={\rm div}\big(\kappa(\FF_{\EPS k},\theta_{\EPS k})\nabla\theta_{\EPS k}\big)
-({\rm div}\,\vv_{\EPS k})\Ent_{\EPS k}
+\TT_1(\FF_{\EPS k},\theta_{\EPS k}){:}\nabla\vv_{\EPS k}
\\&\hspace{7em}\nonumber
+\frac{\nu_1|\ee(\vv_{\EPS k})|^2+\nu_2|\Nabla^2\vv_{\EPS k}|^p}
{1{+}\EPS|\ee(\vv_{\EPS k})|^2{+}\EPS|\Nabla^2\vv_{\EPS k}|^p}
-\frac1k\pdt{\theta_{\EPS k}}
\ \ \ \text{ with }\ \ \Ent_{\EPS k}=\ENT(\FF_{\EPS k},\theta_{\EPS k})\,,
\end{align}
where $\TT_1(\FF,\theta)=\TT(\FF,\theta)-\TT(\FF,0)$  is from 
\eq{stress-split} and $\ENT=\ENT(\FF,\theta)$ is from
\eq{Euler3-enthalpy-finite}.

Relying on  the fact  that $\ENG$ bounded from below due to the
assumption \eq{ass-stress-control} which implies in particular that
$\ENG(\FF,\theta)\ge-1$, the inequality \eq{Euler-engr-finite+} can yield a
first set of a-priori estimates.
To this aim, we are to estimate the right-hand side in \eq{Euler-engr-finite+}.
The issue is
estimation of the gravity force $\varrho\GRAVITY$  when tested by the
velocity $\vv$. For any time instant $t$, this can be estimated by the H\"older/Young
inequality as
\begin{align}\label{Euler-est-of-rhs}
\int_\varOmega\!\varrho_{\EPS k}(t)\GRAVITY(t){\cdot}\vv_{\EPS k}(t)\,\d\xx
&\le\big\|\sqrt{\varrho_{\EPS k}(t)}\big\|_{L^{2}(\varOmega)}
\big\|\sqrt{\varrho_{\EPS k}(t)}\vv_{\EPS k}(t)\big\|_{L^2(\varOmega;\R^d)}^{}
\big\|\GRAVITY(t)\big\|_{L^\infty(\varOmega;\R^d)}^{}
\\&\nonumber\le
\frac12\Big(\big\|\sqrt{\varrho_{\EPS k}(t)}\big\|_{L^2(\varOmega)}^2
\!+\big\|\sqrt{\varrho_{\EPS k}(t)}\vv_{\EPS k}(t)\big\|_{L^2(\varOmega;\R^d)}^2\Big)
\,\big\|\GRAVITY(t)\big\|_{L^\infty(\varOmega;\R^d)}^{}
\\&\!\!\!\stackrel{\scriptsize\eq{mass-conservation}}{=}\!
\big\|\GRAVITY(t)\big\|_{L^\infty(\varOmega;\R^d)}^{}\int_\varOmega\frac{\varrho_0}2
+\frac{\varrho_{\EPS k}(t)}2|\vv_{\EPS k}(t)|^2\,\d\xx\,.
\nonumber\end{align}
The integral on the right-hand side of \eq{Euler-est-of-rhs} can then be
treated by the Gronwall lemma, for which one needs the qualification
of $\GRAVITY$ in \eq{Euler-ass-rhoR}. Here we further use also
the $L^1$-qualification of $h$ in \eq{Euler-ass-theta0}. 
As a result, from \eq{Euler-engr-finite+} we obtain the a-priori estimates
\begin{subequations}\label{Euler-est}\begin{align}
&\big\|\sqrt{\varrho_{\EPS k}}\vv_{\EPS k}\big\|_{L^\infty(I;L^2(\varOmega;\R^d))}^{}\le C\ \ \ \text{ and }\ \ \
\big\|\vv_{\EPS k}\big\|_{L^2(I\times\varOmega;\R^d))}^{}\le C\sqrt k\,,
\label{est-rv2+}
\\&\label{est-phi+}
\big\|\ENG(\FFepsk,\theta_{\EPS k})\big\|_{L^\infty(I;L^1(\varOmega))}^{}\le C\,.
\intertext{Using the energy-controlled-stress condition
\eq{ass-stress-control}, the estimate \eq{est-phi+} implies also}
\label{est-T-L1}
&\big\|\TT(\FFepsk,\theta_{\EPS k})\big\|_{L^\infty(I;L^1(\varOmega;\R^{d\times d}))}^{}\le C\,.
\end{align}\end{subequations}

Now we can use the mechanical-energy balance \eq{Euler-mech-engr-finite+}
together with the symmetry of $\TT$ due to \eq{frame-indifference} and
estimate the right-hand side at each time instant $t$ as
\begin{align}\nonumber
&\int_\varOmega\!\!\!-\TT(\FFepsk,\theta_{\EPS k}){:}\nabla\vv_{\EPS k}\d\xx
=\!\int_\varOmega\!\!\!-\TT(\FFepsk,\theta_{\EPS k}){:}\ee(\vv_{\EPS k})\d\xx
\le\|\TT(\FFepsk,\theta_{\EPS k})\|_{L^1(\varOmega;\R^{d\times d})}^{}\|\ee(\vv_{\EPS k})\|_{L^\infty(\varOmega;\R^{d\times d})}^{}
\\&\hspace{.2em}\nonumber
\le N\|\TT(\FFepsk,\theta_{\EPS k})\|_{L^1(\varOmega;\R^{d\times d})}^{}
\big(\|\ee(\vv_{\EPS k})\|_{L^2(\varOmega;\R^{d\times d})}^{}
+\|\nabla^2\vv_{\EPS k}\|_{L^p(\varOmega;\R^{d\times d\times d})}^{}\big)
\\&\hspace{.2em}\nonumber
\le N\|\TT(\FFepsk,\theta_{\EPS k})\|_{L^1(\varOmega;\R^{d\times d})}^{}
\Big(\frac1{4\delta}+\frac{p-1}{\!\!p^{p/(p-1)}\delta^{1/(p-1)}\!\!}
+\delta\|\ee(\vv_{\EPS k})\|_{L^2(\varOmega;\R^{d\times d})}^2\!
+\delta\|\nabla^2\vv_{\EPS k}\|_{L^p(\varOmega;\R^{d\times d\times d})}^p\Big),
\end{align}
where $N$ denotes the norm of the continuous embedding
$W^{1,p}(\varOmega)\subset L^\infty(\varOmega)$ when the norm
$\|\cdot\|_{L^2(\varOmega)}+\|\nabla\,\cdot\|_{L^p(\varOmega;\R^d)}$ on
$W^{1,p}(\varOmega)$ is used. The constant $\delta>0$ is to be
chosen so small that the last two terms can be absorbed in the left-hand
side of \eq{Euler-mech-engr-finite+}, exploiting also the assumption
\eq{Euler-ass-kappa}. Then \eq{Euler-mech-engr-finite+} gives the additional
estimates
\begin{align}
&\label{est-nabla-v}
\|\ee(\vv_{\EPS k})\|_{L^2(I{\times}\varOmega;\R^{d\times d})}^{}\le C\ \ \ \ \ \text{ and }\ \ \ \ \|\nabla^2\vv_{\EPS k}\|_{L^p(I{\times}\varOmega;\R^{d\times d\times d})}^{}\le C\,.
\end{align}

When assuming $p>d$, the estimate \eq{est-nabla-v} is essential by preventing
the evolution of singularities of the quantities transported by such a smooth
velocity field. The latter estimate in \eq{est-rv2+}, which is not uniform
with $k$, anyhow guarantees that $\vv_{\EPS k}\in L^1(I;W^{2,p}(\varOmega;\R^d))$
for any $k$, which is needed for the proof of the above mentioned properties
of the operators in \eq{transport-equations+}.

We  also have  the following transport-and-evolution equation for $1/\det\FFepsk$ and
a similar equation holds also for the sparsity $1/\varrho_{\EPS k}$, cf.\ \eq{ultimate}
and \eq{cont-eq+}, namely
\begin{align}\label{transport-equations++}
\pdt{}\frac1{\det\FFepsk}=
-{\rm div}\Big(\frac{\vv_{\EPS k}}{\det\FFepsk}\Big)\ \ \ \text{ and }
\ \ \ \pdt{}\frac1{\varrho_{\EPS k}}=\frac{{\rm div}\,\vv_{\EPS k}}{\varrho_{\EPS k}}
-\vv_{\EPS k}{\cdot}\nabla\frac1{\varrho_{\EPS k}}\,.
\end{align}
We can apply the same arguments to the solutions to \eq{transport-equations+}
as we did for \eq{transport-equations} based on \cite[Sect.5]{Roub24TVSE}.
Here, due to the qualification of $\FF_0$ and 
$\varrho_0=\rhoR/\!\det\FF_0$ in \eq{Euler-ass-Fe0}  and \eq{Euler-ass-rhoR}, 
it  follows that  $\varrho_{\EPS k}>0$ and $\det\FFepsk>0$ a.e.\ with
the estimates
\begin{subequations}\label{est+}
\begin{align}
\label{est+Fes}&\|\FFepsk\|_{L^\infty(I;W^{1,r}(\varOmega;\R^{d\times d}))}\le C_r\,,
\ \ \ \:\bigg\|\frac1{\det\FFepsk}\bigg\|_{L^\infty(I;W^{1,r}(\varOmega))}\le C_r\,,
\\&\label{est+rho}\|\varrho_{\EPS k}\|_{L^\infty(I;W^{1,r}(\varOmega))}^{}\le C_r\,,
\ \ \text{ and }\ \ \bigg\|\frac1{\varrho_{\EPS k}}\bigg\|_{L^\infty(I;W^{1,r}(\varOmega))}\!\le C_r
\ \ \ \text{ for any $1\le r<+\infty$};
\intertext{cf.\ also the regularity of the initial conditions that follows
from the assumptions (\ref{Euler-ass}i-k):}
&\nonumber
\nabla\Big(\frac1{\det\FF_0}\Big)=-\,\frac{\det'(\FF_0){:}\nabla\FF_0}
{(\det\FF_0)^2}=-\,\frac{{\rm Cof}\FF_0{:}\nabla\FF_0}{(\det\FF_0)^2}\in L^r(\varOmega;\R^d)\,,
\\&\nonumber
\nabla\varrho_0=\frac{\nabla\rhoR}{\det\FF_0}-\rhoR\frac{{\rm Cof}\FF_0{:}\nabla\FF_0}{(\det\FF_0)^2}\in L^r(\varOmega;\R^d)\,,\ \ \text{ and }\ \ 
\\&\nonumber
\nabla\Big(\frac1{\varrho_0}\Big)=\nabla\Big(\frac{\det\FF_0}{\rhoR}\Big)=
\frac{{\rm Cof}\FF_0{:}\nabla\FF_0}{\rhoR}
-\frac{\det\FF_0\nabla\rhoR}{\varrho_\text{\sc r}^2}\in L^r(\varOmega;\R^d)\,.
\intertext{
From \eq{est+rho} and the former estimate in \eq{est-rv2+}, 
we then have also}
&\big\|\vv_{\EPS k}\big\|_{L^\infty(I;L^2(\varOmega;\R^d))}^{}\le
\big\|\sqrt{\varrho_{\EPS k}}\vv_{\EPS k}\big\|_{L^\infty(I;L^2(\varOmega;\R^d))}^{}
\bigg\|\frac1{\sqrt{\varrho_{\EPS k}}}\bigg\|_{L^\infty(I\times\varOmega)}^{}\le C_r\,,
\label{basic-est-of-v}
\end{align}\end{subequations}
which improves the latter estimate in \eq{est-rv2+}. Realizing the
embedding $W^{1,p}(\varOmega)\subset L^\infty(\varOmega)$ for $p>d$, 
one can complement \eq{est-nabla-v} by
the bound ${\rm ess\,sup}_{t\in I}\|\ee(\vv_{\EPS k}(t))\|_{L^\infty(\varOmega;\R^{d\times d})}^{}\le C_r$.

The estimates \eq{est+Fes}  together with the compact embedding
$W^{1,r}(\varOmega)\subset L^\infty(\varOmega)$ for $r>d$ imply that
both $\FFepsk$ and $1/\det\FFepsk>0$ are bounded on $I{\times}\barOmega$
uniformly in $\EPS$ and $k$. This  also serves for taken a compact set
$K\subset {\rm GL}^+(d)$  which is then  used in the assumptions
(\ref{Euler-ass}e--g). In particular, using \eq{growth-T}, one can complement
also \eq{est-T-L1} by
\begin{align}\label{est-T-Loo}
\bigg\|\frac{\TT(\FFepsk,\theta_{\EPS k})}{1{+}(\theta_{\EPS k}^+)^{1+\ALPHA}}\bigg\|_{L^\infty(I\times\varOmega;\R^{d\times d})}^{}\!\!\le C_{r,K}\,.
\end{align}

It is routine to calculate, when exploiting the continuity equation
\eq{Euler0-thermo-finite+}, that the momentum equation
\eq{Euler1-thermo-finite+} can be written as the balance of the linear
momentum $\varrho\vv$ in the form
\begin{align}\label{momentum-Galerkin}
\pdt{}\big(\varrho\vv\big)={\rm div}\big(\TT(\FF,\theta){+}\Sv
-\varrho\vv{\otimes}\vv\big)+\varrho\GRAVITY-\frac1k\vv\,.
\end{align}
Taking the advantage that the mentioned continuity equation is non-discretized
even in the semi-Galerkin approximation, cf.\ \eq{transport-equations},
we can see that the approximate solution solves also \eq{momentum-Galerkin}
in its Galerkin approximation with the same finite-dimensional space $V_k$
as used for the Galerkin approximation of \eq{Euler1-thermo-finite+}.

\medskip\noindent{\it Step 3: existence of semi-Galerkin solution and an $L^2$-estimate of $\nabla\theta_{\EPS k}$}.
We still test \eq{Euler3-thermo-finite+} by $\theta_{\EPS k}$ which gives an
estimate for $\nabla\theta_{\EPS k}$ because of the regularization of the
right-hand side in \eq{Euler3-thermo-finite+}. We have the pointwise estimate
$|\TT(\FFepsk,\theta_{\EPS k})|\le C_{r,K}(1{+}(\theta_{\EPS k}^+)^{1+\ALPHA})$
with $C_{r,K}$ from \eq{est-T-Loo} and the same type of estimate holds also
for $\TT_1(\FF,\theta)=\TT(\FF,\theta)-\TT(\FF,0)$.
Therefore, exploiting the embedding $W^{1,p}(\varOmega)\subset L^\infty(\varOmega)$
and the symmetry of the Cauchy stress $\TT_1$ granted by \eq{frame-indifference},
we have also at each time instant $t$ the pointwise  estimate for the power of
the conservative part of the Cauchy stress (i.e.\ the internal-energy source) as
\begin{align}
\big|\TT_1(\FFepsk(t),\theta_{\EPS k}(t)){:}\nabla\vv_{\EPS k}(t)\big|
=\big|\TT_1(\FFepsk(t),\theta_{\EPS k}(t)){:}\ee(\vv_{\EPS k}(t))\big|
\le C(t)\big(1{+}(\theta_{\EPS k}^+(t))^{1+\ALPHA}\big)
\label{est-adibatic-heat}\end{align}
with some $C\in L^p(I)$ depending on the estimates 
\eq{est-T-Loo}.

Let us consider $\wt\ENG(\FF,\theta)$ so that
$\wt\ENG_\theta'(\FF,\theta)=\theta c(\FF,\theta)$ with the heat capacity
$c(\FF,\theta)=\ENG_\theta'(\FF,\theta)$. Specifically, we put
$\wt\ENG(\FF,\theta):=\int_0^\theta\wh\theta c(\FF,\wh\theta)\,\d\wh\theta$ or,
equivalently,
\begin{align}\label{primitive-theta-c}
\wt\ENG(\FF,\theta)=\int_0^1\!\!r\theta^2c(\FF,r\theta)\,\d r\,.
\end{align}
Let us further put $\wt\ENT(\FF,\theta):=\wt\ENG(\FF,\theta)-\theta\,\ENT(\FF,\theta)$.
From \eq{growth-int-eng}, we can see that $|\wt\ENT(\FF,\theta)|\le2
C_K^{}(1{+}(\theta^+)^{2+\ALPHA})$ with  $C_K^{}$ referring to 
the compact subset  $K$ 
of ${\rm GL}^+(d)$ respecting the already obtained estimates \eq{est+Fes}.
Then, like in \cite{Roub24TVSE}, we have the calculus
\begin{align*}
\theta_{\EPS k}\big(\ENT(\FFepsk,\theta_{\EPS k})\big)\!\!\DTk{^{^{}}}\ 
&=\theta_{\EPS k}\ENT_\theta'(\FFepsk,\theta_{\EPS k})\DTk\theta_{\EPS k}
+\theta_{\EPS k}\ENT_\FF'(\FFepsk,\theta_{\EPS k}){:}\!\!\DTk{\,\,\FFepsk}
\\[-.3em]&=\wt\ENG_\theta'(\FFepsk,\theta_{\EPS k})\DTk\theta_{\EPS k}
+\theta_{\EPS k}\ENT_\FF'(\FFepsk,\theta_{\EPS k}){:}\!\!\DTk{\,\,\FFepsk}
\\[-.3em]&=\big({\wt\ENG}(\FFepsk,\theta_{\EPS k})\big)\!\!\DTk{^{^{}}}\, 
-\wt\ENG_\FF'(\FFepsk,\theta_{\EPS k}){:}\!\!\DTk{\,\,\FFepsk}
+\theta_{\EPS k}\ENT_\FF'(\FFepsk,\theta_{\EPS k}){:}\!\!\DTk{\,\,\FFepsk}
\\[-.3em]&\!\!\!\stackrel{\scriptsize\eq{transport-equations}}{=}\!\!
\big({\wt\ENG}(\FFepsk,\theta_{\EPS k})\big)\!\!\DTk{^{^{}}}\, 
-\big(\wt\ENG_\FF'(\FFepsk,\theta_{\EPS k})
-\theta_{\EPS k}\ENT_\FF'(\FFepsk,\theta_{\EPS k})\big){:}\big((\nabla\vv_{\EPS k})\FFepsk\big)
\\[-.3em]&=
\big({\wt\ENG}(\FFepsk,\theta_{\EPS k})\big)\!\!\DTk{^{^{}}}\, 
-\wt\ENT_\FF'(\FFepsk,\theta_{\EPS k})\FFepsk^\top{:}\nabla\vv_{\EPS k}\,,
\end{align*}
{}where the notation $(\cdot)\!\DTk{^{}\,}\,$ is from  \eq{def-of-DTek}. 
In terms of $\wt\ENG$ and $\wt\ENT$, we can write the Galerkin approximation of
\eq{Euler3-enthalpy-finite-reg} tested by $\theta_{\EPS k}$ as
\begin{align}\label{heat-test-theta}
&\int_\varOmega\wt\ENG(\FFepsk(t),\theta_{\EPS k}(t))+\frac1k|\theta_{\EPS k}(t)|^2\,\d\xx
+\!\int_0^t\!\!\int_\varOmega\kappa(\FFepsk,\theta_{\EPS k})|\nabla\theta_{\EPS k}|^2\,\d\xx
\d t
\\[-.1em]&\hspace{.5em}\nonumber
=\int_\varOmega\wt\ENG(\FF_0,\theta_{0\EPS})\,\d\xx
+\!\int_0^t\!\!\int_\varOmega\bigg(\Big(\frac{\nu_1|\ee(\vv_{\EPS k})|^2
+\nu_2|\nabla^2\vv_{\EPS k}|^p\!\!}
{1{+}\EPS|\ee(\vv_{\EPS k})|^2{+}\EPS|\nabla^2\vv_{\EPS k}|^p}
+\TT_1(\FFepsk,\theta_{\EPS k}){:}\nabla\vv_{\EPS k}\Big)\theta_{\EPS k}
\\[-.2em]&\hspace{4.8em}\nonumber
+\wt\ENT_\FF'(\FFepsk,\theta_{\EPS k})\FFepsk^\top{:}\nabla\vv_{\EPS k}
+\wt\ENT(\FFepsk,\theta_{\EPS k}){\rm div}\,\vv_{\EPS k}\!\bigg)\,\d\xx\d t
+\!\int_0^t\!\!\int_\varGamma h_\EPS\theta_{\EPS k}\,\d S\d t
\\[-.1em]&\hspace{.5em}\nonumber
\le\int_\varOmega\wt\ENG(\FF_0,\theta_{0\EPS})\,\d\xx
+\!\int_0^t\!\!\int_\varOmega\bigg(\frac{\nu_1|\ee(\vv_{\EPS k})|^2
+\nu_2|\nabla^2\vv_{\EPS k}|^p\!\!}
{1{+}\EPS|\ee(\vv_{\EPS k})|^2{+}\EPS|\nabla^2\vv_{\EPS k}|^p}\:\theta_{\EPS k}^+
+\big(\TT_1(\FFepsk,\theta_{\EPS k}){:}\nabla\vv_{\EPS k}\big)\theta_{\EPS k}^+
\\[-.5em]&\hspace{4.8em}
+\wt\ENT_\FF'(\FFepsk,\theta_{\EPS k})\FFepsk^\top{:}\nabla\vv_{\EPS k}
+\wt\ENT(\FFepsk,\theta_{\EPS k}){\rm div}\,\vv_{\EPS k}\!\bigg)\,\d\xx\d t
+\!\int_0^t\!\!\int_\varGamma 
\frac{h\,\theta_{\EPS k}^+}{1{+}\EPS h}\,\d S\d t\,.
\nonumber
\end{align}
We also used that $h\ge0$ as assumed in \eq{Euler-ass-theta0} and that
$\TT_1(\cdot,\theta)=\bm0$ for $\theta\le0$, as noted in Remark~\ref{rem-split}.
Here we use the regularization in \eq{Euler-thermodynam-IC-reg} which makes
$\wt\ENG(\FF_0,\theta_{0\EPS})$ integrable.
From \eq{est-adibatic-heat}, at each time instant, we can read the estimate
$|\TT_1(\FFepsk,\theta_{\EPS k}){:}(\nabla\vv_{\EPS k})\theta_{\EPS k}|
\le C(t)(1{+}(\theta_{\EPS k}^+)^{2+\ALPHA})$ with $C\in L^p(I)$.
Since $\wt\ENG_\FF'(\FF,\theta)$ is a component-wise primitive function
(antiderivative) of $\theta\mapsto\theta c_\FF'(\FF,\theta)=\theta\ENG_{\FF\theta}''(\FF,\theta)$,
cf.\ \eq{primitive-theta-c}, we can rely on the assumptions \eq{Euler-ass-psi-1+++}
and therefore $|\wt\ENG_{\FF}'(\FF,\theta)|\le 
C_K^{}(1{+}(\theta^+)^{2+\ALPHA})$ and thus also $|\wt\ENT_\FF'(\FF,\theta)|\le
2C_K^{}(1{+}(\theta^+)^{2+\ALPHA})$. The coercivity 
\eq{Euler-ass-psi-1} yields the coercivity
$\wt\ENG(\FF,\theta)\ge c_K^{}(\theta^+)^{2+\ALPHA}/(2{+}\ALPHA)$, 
which allows for usage of the Gronwall inequality for \eq{heat-test-theta}.
The estimation of the boundary term in \eq{heat-test-theta}
exploits $\int_\Gamma h\theta_{\EPS k}^+/(1{+}\EPS h)\,\d S\le\EPS^{-1}
\|\theta_{\EPS k}^+\|_{L^1(\varGamma)}\le
C_{\ALPHA,\EPS}(1{+}\|\theta_{\EPS k}^+\|_{L^{2+\ALPHA}(\varOmega)}^{2+\ALPHA})
+\frac12c_K\|\nabla\theta_{\EPS k}\|_{L^2(\varOmega;\R^d)}^2$ with $c_K$ referring to
\eq{Euler-ass-kappa} and with some $C_{\ALPHA,\EPS}$ sufficiently large, so that the
mentioned coercivity of $\wt\ENG(\FF,\cdot)$ on $\R^+$ and the Gronwall inequality
can be used.
From the already obtained estimates for the Galerkin approximation, we can thus
read some a-priori estimated also for $\theta_{\EPS k}$, namely
\begin{align}\label{est-theta+}
&\big\|\nabla\theta_{\EPS k}\big\|_{L^2(I\times\varOmega;\R^d)}^{}\le C, \ \ \
\big\|\theta_{\EPS k}^+\big\|_{L^\infty(I;L^{2+\ALPHA}(\varOmega))\,\cap\,L^2(I;H^1(\varOmega))}^{}\le C,
\ \text{ and }\ \big\|\theta_{\EPS k}\big\|_{L^\infty(I;L^2(\varOmega))}^{}\le C\sqrt{k}\,.
\end{align}
Note that, disregarding $\nabla\theta^-$, we do not have any uniform control on the
negative part of temperature which may be indeed nonvanishing in the Galerkin
approximation, but it does not harm the arguments in the proof because the
data $\TT(\FF,\cdot)$ and $\kappa(\FF,\cdot)$ are insensitive to possible negative
values of temperature, cf.\ the assumption \eq{Euler-ass-psi-increas}. Anyhow, the
last $L^\infty$-estimate of $\theta_{\EPS k}$ in \eq{est-theta+}, although
being nonuniform in $k$, together with the other $L^\infty$-estimates
(\ref{est+}a--c) allow us to apply the successive-prolongation arguments and to show
the existence of a semi-discrete solution
$(\varrho_{\EPS k},\vv_{\EPS k},\FFepsk,\theta_{\EPS k})$ not only locally
in time but globally on the whole time interval $I$.

Furthermore, from \eq{est-theta+}, we can read an estimate for
\begin{align}\nonumber
\nabla\Ent_{\EPS k}&=\ENT_{\FF}'(\FFepsk,\theta_{\EPS k}){:}\nabla\FFepsk+
\ENT_{\theta}'(\FFepsk,\theta_{\EPS k})\nabla\theta_{\EPS k}\\&=
\big(\psi_{\FF}'(\FFepsk,\theta_{\EPS k})-\theta\psi_{\FF\theta}''(\FFepsk,\theta_{\EPS k})
-\psi_{\FF}'(\FFepsk,0)\big){:}\nabla\FFepsk+c(\FFepsk,\theta_{\EPS k})\nabla\theta_{\EPS k}
\nonumber\end{align}
with the heat capacity
$c(\FFepsk,\theta_{\EPS k})=-\theta_{\EPS k}\psi_{\theta\theta}''(\FFepsk,\theta_{\EPS k})$.
Specifically, due to \eq{Euler-ass-psi-1+++} with \eq{Euler-ass-psi-increas},
\eq{est+Fes}, and \eq{est-theta+}, we have $\ENT_{\FF}'(\FFepsk,\theta_{\EPS k})$
bounded in $L^{(2+\ALPHA)/(1+\ALPHA)}(I{\times}\varOmega;\R^{d\times d})$
while $\nabla\FFepsk$ is bounded in $L^{r}(I{\times}\varOmega;\R^{d\times d\times d})$
for any (arbitrarily big) $1\le r<+\infty$. Thus
$\ENG_{\FF}'(\FFepsk,\theta_{\EPS k}){:}\nabla\FFepsk$ is
bounded in $L^{r(2+\ALPHA)/(2+\ALPHA+r+r\ALPHA)}(I{\times}\varOmega;\R^d)$.
Realizing that $c(\FFepsk,\theta_{\EPS k})\le C(1{+}(\theta^+)^\ALPHA)$ as actually
assumed by the first condition in \eq{Euler-ass-psi-1+++}, from \eq{est-theta+}
we have $c(\FFepsk,\theta_{\EPS k})\nabla\theta_{\EPS k}$ bounded in
$L^{(4+2\ALPHA)/(2+3\ALPHA)}(I{\times}\varOmega;\R^d)$, so that altogether (considering
$r$ sufficiently big) we have 
\begin{subequations}
\begin{align}\label{est-w+}
&\big\|\nabla\Ent_{\EPS k}\big\|_{L^{(4+2\ALPHA)/(2+3\ALPHA)}(I\times\varOmega;\R^d)}^{}\le C
\intertext{and, due to the first estimate in \eq{est-theta+}, we have also}
&\Big\|\nabla\Big(\Ent_{\EPS k}+\frac{\theta_{\EPS k}}k\Big)\Big\|_{L^{(4+2\ALPHA)/(2+3\ALPHA)}(I\times\varOmega;\R^d)}^{}\le C
\,.
\end{align}
\end{subequations}

\medskip\noindent{\it Step 4: Convergence for $k\to\infty$}.
Using the Banach selection principle, we can extract some subsequence of
$\{(\varrho_{\EPS k},\vv_{\EPS k},\FFepsk,\theta_{\EPS k}^+,\Ent_{\EPS k}{+}\theta_{\EPS k}/k)\}_{k\in\N}$
and its limit
$(\varrho_{\EPS},\vv_{\EPS},\FF_{\EPS},\theta_{\EPS},\Ent_{\EPS}):I\to W^{1,r}(\varOmega)\times
L^2(\varOmega;\R^d)\times W^{1,r}(\varOmega;\R^{d\times d})\times
L^{2+\ALPHA}(\varOmega)\times L^1(\varOmega)$
such that, for any $1\le r<+\infty$, we have
\begin{subequations}\label{Euler-disc-sln}
\begin{align}
&\!\!\varrho_{\EPS k}\to\varrho_{\EPS}&&\text{weakly* in $\
L^\infty(I;W^{1,r}(\varOmega))\,\cap\,W^{1,
p}(I;L^r(\varOmega))$}\,,
\\\label{Euler-disc-sln-v}
&\!\!\vv_{\EPS k}\to\vv_{\EPS}&&\text{weakly* in $\
L^\infty(I;L^2(\varOmega;\R^d))\cap
L^{
p}(I;W^{2,p}(\varOmega;\R^d))$,}\!\!&&
\\\label{Euler-disc-sln-F}
&\!\!\FFepsk\to\FF_{\EPS}
\!\!\!&&\text{weakly* in $\ L^\infty(I;W^{1,r}(\varOmega;\R^{d\times d}))\,\cap\,
W^{1,
p}(I;L^r(\varOmega;\R^{d\times d}))$,}\!\!
\\\label{Euler-disc-sln-theta}
&\!\!\theta_{\EPS k}^+\to\,\theta_{\EPS}\!\!\!&&\text{weakly* in $\
L^\infty(I;L^{2+\ALPHA}(\varOmega))\,\cap\,L^2(I;H^1(\varOmega))$,}\!\!
\\[-.4em]
&\!\!\Ent_{\EPS k}\,{+}\,\frac{\theta_{\EPS k}}k\to\Ent_{\EPS}\hspace{-1.5em}&&
\text{weakly in $\ L^{(4+2\ALPHA)/(2+3\ALPHA)}(I;
W^{(4+2\ALPHA)/(2+3\ALPHA)}(\varOmega))$}\,.
\label{Euler-disc-sln-w++}
\end{align}\end{subequations}
In \eq{Euler-disc-sln-F}, we used also the estimates \eq{est-nabla-v} and
(\ref{est+}a,c), which yields the bound of
$\pdt{}\FFepsk=(\nabla\vv_{\EPS k})\FFepsk-(\vv_{\EPS k}{\cdot}\nabla)\FFepsk$ in
$L^p(I;L^r(\varOmega;\R^{d\times d}))$.  Similarly, using (\ref{est+}b,c), we
obtain also the bound for
$\pdt{}\varrho_{\EPS k}=({\rm div}\vv_{\EPS k})\varrho_{\EPS k}-\vv_{\EPS k}{\cdot}\nabla\varrho_{\EPS k}$ in $L^p(I;L^r(\varOmega))$.  
By the Aubin-Lions lemma here considering $r>d$, we also have that
\begin{align}\label{rho-conv}
&\varrho_{\EPS k}\to\varrho_{\EPS}\ \text{ strongly in }\ C(I{\times}\barOmega) 
\ \ \text{ and }\ \ \FFepsk\to\FF_{\EPS}\ \text{ strongly in
$C(I{\times}\barOmega;\R^{d\times d})$\,.}
\end{align}
This already allows for the limit passage in the evolution equations
\eq{transport-equations}.

Further, by comparison in the equation \eq{Euler3-enthalpy-finite-reg} with the
boundary condition \eq{Euler-thermodynam-BC-2-reg} in its Galerkin
approximation, we obtain a bound  for  $\pdt{}\Eng_{\EPS k}$ in the seminorms
$|\cdot|_l$ on $L^2(I;H^1(\varOmega)^*)$ arising from this Galerkin
approximation:
\begin{align}\nonumber
\big|f\big|_l^{}:=
\sup\limits_{\|\widetilde\theta\|_{L^2(I;H^1(\varOmega))}^{}\le1\,,\ \widetilde\theta(t)\in Z_l\ \text{for }t\in I}
\int_0^T\!\!\!\int_\varOmega f\widetilde\theta\,\d\xx\d t\,.
\end{align}
Specifically, for any $k\ge l$, we can estimate
\begin{align}\nonumber
&\bigg|\pdt{}\Big(\Ent_{\EPS k}{+}\frac{\theta_{\EPS k}}k\Big)\bigg|_l^{}=\!\!\!\!\!
\sup\limits_{\begin{array}{c}{\scriptstyle{\widetilde\theta(t)\in Z_l\ \text{for }t\in I}}\\[-.3em]
{\scriptstyle{\|\widetilde\theta\|_{L^2(I;H^1(\varOmega))}^{}\le1}}
\end{array}
}\!\!\!\!
\int_0^T\!\!\!\int_\varOmega\!\bigg(\Big(\frac{\nu_1|\ee(\vv_{\EPS k})|^2\!
+\nu_2|\nabla^2\vv_{\EPS k}|^p}{1{+}\EPS|\ee(\vv_{\EPS k})|^2{+}\EPS|\nabla^2\vv_{\EPS k}|^p}
+\TT_1(\FFepsk,\theta_{\EPS k}){:}\nabla\vv_{\EPS k}
\Big)\widetilde\theta
\\[-2em]&\hspace{14.5em}\nonumber
+\big(\Ent_{\EPS k}\vv_{\EPS k}{-}\kappa(\FFepsk,\theta_{\EPS k})\nabla\theta_{\EPS k}\big)
{\cdot}\nabla\widetilde\theta\,\bigg)\d\xx\d t
+\!\int_0^T\!\!\!\int_\varGamma\!
h_\EPS\widetilde\theta\,\d S\d t\le C
\end{align}
with some $C$ depending on the estimates  \eq{est-T-Loo} and \eq{est-theta+} but
independent  of  $l\in \N$. Thus, by \eq{Euler-disc-sln-w++} and by a generalized
Aubin-Lions theorem \cite[Ch.8]{Roub13NPDE}, we obtain
\begin{align}
&&&\label{u-conv}
\Ent_{\EPS k}\,{+}\,\frac{\theta_{\EPS k}}k\to \Ent_{\EPS}\hspace*{-0em}&&
\hspace*{-1em}\text{strongly in $L^{(4+2\ALPHA)/(2+3\ALPHA)}(I{\times}\varOmega)$};&&
\end{align}
actually, the exponent in  \eq{u-conv}  taken from
\eq{Euler-disc-sln-w++} could be improved when interpolating also the estimate
\eq{est-phi+}. Due to the last estimate in \eq{est-theta+},
$\|\theta_{\EPS k}/k\|_{L^2(I{\times}\varOmega)}^{}=\mathscr{O}(1/\sqrt{k})\to0$, and
therefore we have also $\Ent_{\EPS k}\to \Ent_{\EPS}$
strongly in $L^{(4+2\ALPHA)/(2+3\ALPHA)}(I{\times}\varOmega)$. Moreover, in view
of the continuity of $\ENG(\cdot,0)$ on the compact $K$  in which 
$\FF_{\EPS k}$ are valued, \eq{rho-conv} implies also
\begin{align}
&&&\label{w-conv}
\Eng_{\EPS k}=\Ent_{\EPS k}{-}\,\ENG(\FF_{\EPS k},0)\,\to\,\Ent_{\EPS}{-}\,\ENG(\FF_{\EPS},0)
=:\Eng_{\EPS}\hspace*{-0em}&&\text{strongly in $L^{(4+2\ALPHA)/(2+3\ALPHA)}(I{\times}\varOmega)$}.&&
\end{align}

We unfortunately cannot read convergence of
$\theta_{\EPS k}\in[\ENT(\FFepsk,\cdot)]^{-1}(\Ent_{\EPS k})$ from the mentioned
convergence of $\Ent_{\EPS k}$ because the inverse $[\ENT(\FF,\cdot)]^{-1}$
of $\ENT(\FF,\cdot):\R\to\R$ is not continuous, being set-valued at 0, specifically
$[\ENT(\FF,\cdot)]^{-1}(0)=(-\infty,0]$ while
$[\ENT(\FF,\cdot)]^{-1}(\Ent)=\emptyset$ for $\Ent<0$. Anyhow, 
since $\ENG(\FF,\cdot):\R\to\R^+$ is assumed increasing and coercive
on $\R^+$, cf.\ \eq{Euler-ass-psi-increas} and \eq{Euler-ass-psi-1},
so is $\ENT(\FF,\cdot)=\ENG(\FF,\cdot)-\ENG(\FF,0)$, and thus
there exists a continuous inverse $[\ENT(\FF,\cdot)]^{-1}:(0
,+\infty)\to\R$. This analytical peculiarity can be circumvented by
defining a modification for $\Ent<0$ as a continuous extension to
obtain the single-valued continuous function
\begin{align}
&\big[\ENT(\FF,\cdot)\big]_{\text{\sc m}}^{-1}:\R\to\R^+:
\Ent\mapsto\begin{cases}[\ENT(\FF,\cdot)]^{-1}(\Ent)\!\!
&\text{for }\ \Ent>0
\,,\\\qquad0&\text{for }\ \Ent\le0
\,;\end{cases}
\label{inverse-ext}\end{align}
cf.\ Figure~\ref{various-nonlinearities}-right for illustration.
\begin{figure}[ht]
\begin{center}
\psfrag{t}{\scriptsize $\theta$}
\psfrag{w1}{\scriptsize $\Ent=\ENT(\FF,\theta)$}
\psfrag{w2}{\scriptsize $\wt\ENG(\FF,\theta)$}
\psfrag{w-inv}{\scriptsize $\big[\ENT(\FF,\cdot)\big]^{-1}$}
\psfrag{w-ext}{\scriptsize $\big[\ENT(\FF,\cdot)\big]_{\text{\sc m}}^{-1}$}
\psfrag{w}{\scriptsize $\Ent$}
\psfrag{c}{\scriptsize $c(\FF,\theta)$}
\includegraphics[width=.98\textwidth]{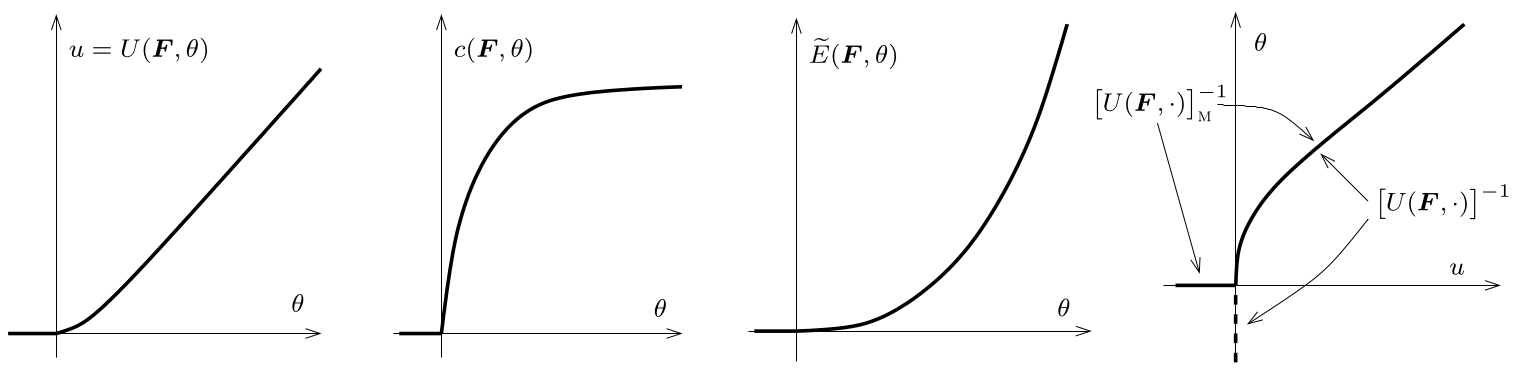}
\end{center}
\vspace*{-.6em}
\caption{
{\sl A schematic illustration of continuous extensions for negative
temperatures of the thermal part of the internal energy $\Ent=\ENT(\FF,\theta)$,
the heat capacity $c=c(\FF,\theta)$, the nonlinearity $\wt\ENG$ from
\eq{primitive-theta-c}, and the modified continuous inverse
$\big[\ENT(\FF,\cdot)\big]_{\text{\sc m}}^{-1}$ of $\ENT(\FF,\cdot)$ as
defined by \eq{inverse-ext}; here $\FF$ is considered fixed.
}
}
\label{various-nonlinearities}
\end{figure}

Since surely $\Ent_{\EPS k}\ge0$, we can write
$\theta_{\EPS k}^+=[\ENT(\FFepsk,\cdot)]_{\text{\sc m}}^{-1}(\Ent_{\EPS k})$. Thanks to the
continuity of
$(\FF,\Ent)\mapsto[\ENT(\FF,\cdot)]_{\text{\sc m}}^{-1}(\Ent):\R^{d\times d}\times\R^+\to\R$
 derived by  the inverse-function theorem  when  using the assumption
$\ENT_\theta'=\ENG_\theta'>0$ in \eq{Euler-ass-psi-1}, we can improve
\eq{Euler-disc-sln-theta}  by 
\begin{align}
&&&\label{z-conv}
\theta_{\EPS k}^+\to \theta_{\EPS}=\big[\ENT(\FF_{\EPS},\cdot)\big]_{\text{\sc m}}^{-1}(\Ent_{\EPS})
&&\text{strongly in $L^{2+4/d}(I{\times}\varOmega)$};&&&&&&
\end{align}
here we have used the Gagliardo-Nirenberg interpolation for the embedding
$L^\infty(I;L^2(\varOmega))\cap L^2(I;H^1(\varOmega))\subset L^{2+4/d}(I{\times}\varOmega)$.
Let us emphasize that we do not have any direct information about
$\pdt{}\theta_{\EPS k}$ neither for $\pdt{}\theta_{\EPS k}^+$ so that we could not use
the Aubin-Lions arguments directly for $\{\theta_{\EPS k}\}_{k\in\N}$ neither for
$\{\theta_{\EPS k}^+\}_{k\in\N}$.

Thus, by the continuity of the corresponding Nemytski\u{\i} (or here simply
superposition) mappings, also the  conservative part of the  Cauchy stress $\TT$
as well as the coefficient $\kappa$ converge strongly, namely
\begin{subequations}\label{Euler-weak+}
\begin{align}
&
\TT(\FFepsk,\theta_{\EPS k})\to \TT(\FF_{\EPS},\theta_{\EPS})
\hspace*{-2em}&&\hspace*{0em}\text{strongly in
$L^r(I{\times}\varOmega;\R_{\rm sym}^{d\times d}),\ \ 1\le r<
\smaller{\frac{2d+4}{d(1{+}\ALPHA)}}$, and}
\label{Euler-T-strong-conv}
\\[-.4em]&\kappa(\FFepsk,\theta_{\EPS k})\to\kappa(\FF_{\EPS},\theta_{\EPS})\hspace*{-2em}&&
\hspace*{0em}\text{strongly in $L^r(I{\times}\varOmega),\ \ 1\le r<\infty$.}\label{Euler-weak+kappa}
\intertext{Here the assumption \eq{Euler-ass-psi-increas} has been employed for
$\TT(\FFepsk,\theta_{\EPS k})=\TT(\FFepsk,\theta_{\EPS k}^+)$ and 
$\kappa(\FFepsk,\theta_{\EPS k})=\kappa(\FFepsk,\theta_{\EPS k}^+)$.
 Moreover,}
&\label{rho.v-conv}
\varrho_{\EPS k}\vv_{\EPS k}\to\varrho_\EPS\vv_\EPS&&
\hspace*{0em}\text{strongly in $L^a(I{\times}\varOmega;\R^d)$ with any $1\le a<4$}
\intertext{{}due to the (generalized) Aubin-Lions compacteness
theorem, relying on that $\nabla(\varrho_{\EPS k}\vv_{\EPS k})
=\varrho_{\EPS k}\nabla\vv_{\EPS k}+\nabla\varrho_{\EPS k}{\otimes}\vv_{\EPS k}$
is bounded in $L^\infty(I;L^r(\varOmega;\R^{d\times d}))$ due to \eq{est-nabla-v}
and (\ref{est+}b,c) and on that, exploiting \eq{momentum-Galerkin}, we have
an information about $\pdt{}(\varrho_{\EPS k}\vv_{\EPS k})$
in the seminorms induced by the Galerkin approximation in the
space $L^2(I;H^1(\varOmega;\R^d)^*)+L^{p'}(I;W^{2,p}(\Omega;\R^d)^*)$. From this,
for the approximate velocity
$\vv_{\EPS k}=(1/\varrho_{\EPS k})\varrho_{\EPS k}\vv_{\EPS k}$ itself, we also have}
&\label{v-conv}
\vv_{\EPS k}\to\vv_\EPS&&
\hspace*{0em}\text{strongly in $L^a(I{\times}\varOmega;\R^d)$ with any $1\le a<4$}
\intertext{{}when using \eq{rho-conv} with \eq{est+rho} for the strong
convergence of $1/\varrho_{\EPS k}$. Cf.\ also
\cite[Formulas (3.35)-(3.36)]{Roub24TVSE}.  For the limit passage in the
dissipative heat power, we will  also need } 
&\ee(\vv_{\EPS k})\to\ee(\vv_{\EPS})\hspace*{0em}&&\text{strongly in $\ L^2(I{\times}\varOmega;\R^{d\times d})\ $ and}
\\&
\Nabla^2\vv_{\EPS k}\to\Nabla^2\vv_{\EPS}\hspace*{0em}&&\text{strongly in $\ L^p(I{\times}\varOmega;\R^{d\times d\times d})\,$,}
\end{align}\end{subequations}
which is based on the uniform monotonicity of the dissipative stress, i.e., of
the quasilinear operator
$\vv\mapsto{\rm div}({\rm div}(\nu_2|\Nabla^2\vv|^{p-2}\Nabla^2\vv)-\nu_1\ee(\vv)))$.
Thus, to prove (\ref{Euler-weak+}e,f), 
we use the Galerkin approximation of the regularized momentum equation
\eq{Euler1-weak-Galerkin} tested by $\widetilde\vv=\vv_{\EPS k}-\widetilde\vv_k$
with $\widetilde\vv_k:I\to V_k$ an approximation of $\vv_\EPS$ in the sense that
$\widetilde\vv_k\to\vv_\EPS$ strongly in $L^p(I;W^{2,p}(\varOmega;\R^d))$
and $\widetilde\vv_k(0)\to\vv_0$ strongly in $L^2(\varOmega;\R^d)$
for $k\to\infty$ while $\widetilde\vv_k{\cdot}\nn=0$ on $I{\times}\varGamma$
and $\{\pdt{}\widetilde\vv_k\}_{k\in\N}$ is bounded in
$L^{p'}\!(I;W^{2,p}(\varOmega;\R^d)^*)$. Using the perturbed continuity equation
$\pdt{}\varrho_{\EPS k}+{\rm div}(\varrho_{\EPS k}(\vv_{\EPS k}{-}\widetilde\vv_k))=
-{\rm div}(\varrho_{\EPS k}\widetilde\vv_k)$ tested by $|\vv_{\EPS k}{-}\widetilde\vv_k|^2/2$ to derive the identity
\begin{align}\nonumber
\int_0^T\!\!\!\int_\varOmega&
\varrho_{\EPS k}\DTk\vv_{\EPS k}
{\cdot}(\vv_{\EPS k}{-}\widetilde\vv_k)\,\d \xx\d t
=\int_0^T\!\!\!\int_\varOmega\bigg(\varrho_{\EPS k}\Big(\pdt{}
(\vv_{\EPS k}{-}\widetilde\vv_k)
{\cdot}(\vv_{\EPS k}{-}\widetilde\vv_k)+
((\vv_{\EPS k}{-}\widetilde\vv_k){\cdot}\nabla)
(\vv_{\EPS k}{-}\widetilde\vv_k)\Big)
\\[-.5em]&\hspace{15em}\nonumber
+\varrho_{\EPS k}\Big(\pdt{\widetilde\vv_k}{\cdot}(\vv_{\EPS k}{-}\widetilde\vv_k)+
(\widetilde\vv_k{\cdot}\nabla)
(\vv_{\EPS k}{-}\widetilde\vv_k)\Big)\bigg)\,\d \xx\d t
\\[-.6em]&\nonumber
=\int_0^T\!\!\!\int_\varOmega\varrho_{\EPS k}\Big(\pdt{\widetilde\vv_k}{\cdot}(\vv_{\EPS k}{-}\widetilde\vv_k)+
(\widetilde\vv_k{\cdot}\nabla)
(\vv_{\EPS k}{-}\widetilde\vv_k)\Big)
+\frac12{\rm div}(\varrho_{\EPS k}\widetilde\vv_k)|\vv_{\EPS k}{-}\widetilde\vv_k|^2\,\d \xx\d t
\\[-.3em]\nonumber&\hspace{15em}
+\int_\varOmega\frac{\varrho_{\EPS k}(T)}2|\vv_{\EPS k}(T){-}\widetilde\vv_k(T)|^2
-\frac{\varrho_0}2|\vv_0{-}\widetilde\vv_k(0)|^2\,\d \xx\,,
\end{align}
we can estimate
\begin{align}\label{strong-hyper+}
&\!\!\!\nu_1\|\ee(\vv_{\EPS k}{-}\vv_\EPS)\|_{L^2(I{\times}\varOmega;\R^{d\times d})}^2
+\widetilde\nu_2\|\nabla^2(\vv_{\EPS k}{-}\vv_\EPS)\|_{L^p(I{\times}\varOmega;\R^{d\times d\times d})}^p
 \\[-.1em]&\nonumber
 \le\int_\varOmega\!\frac{\varrho_{\EPS k}(T)}2\big|\vv_{\EPS k}(T){-}\widetilde\vv_k(T)\big|^2\,\d\xx
+\int_0^T\!\!\!\int_\varOmega\!\bigg(\nu_1\big|\ee(\vv_{\EPS k}{-}\vv_\EPS)\big|^2
\\[-.4em]&\nonumber\hspace{6em}
 +\nu_2\big(|\nabla^2\vv_{\EPS k}|^{p-2}\nabla^2\vv_{\EPS k}
-|\nabla^2\vv_\EPS|^{p-2}\nabla^2\vv_\EPS\big)\Vdots
  \nabla^2(\vv_{\EPS k}{-}\vv_\EPS)\bigg)\,\d\xx\d t
 \\[-.4em]&=\nonumber
 \int_0^T\!\!\!\int_\varOmega\bigg(\varrho_{\EPS k}\GRAVITY{\cdot}(\vv_{\EPS k}{-}\widetilde\vv_k)
-\TT(\FFepsk,\theta_{\EPS k}){:}\ee(\vv_{\EPS k}{-}\widetilde\vv_k)
-\frac{\vv_{\EPS k}}k{\cdot}(\vv_{\EPS k}{-}\widetilde\vv_k)
\\[.1em]&\nonumber\hspace{2em}
-\nu_1\ee(\vv_\EPS){:}\ee(\vv_{\EPS k}{-}\vv_\EPS)
+\nu_1\ee(\vv_{\EPS k}){:}\ee(\widetilde\vv_k{-}\vv_\EPS)
-\nu_2|\nabla^2\vv_\EPS|^{p-2}\nabla^2\vv_\EPS\Vdots
 \nabla^2(\vv_{\EPS k}{-}\vv_\EPS)
\\[-.1em]&\nonumber\hspace*{2em}
+\nu_2|\nabla^2\vv_{\EPS k}|^{p-2}\nabla^2\vv_{\EPS k}\Vdots
  \nabla^2(\widetilde\vv_k{-}\vv_\EPS)
-\varrho_{\EPS k}\Big(\pdt{\widetilde\vv_k}{\cdot}(\vv_{\EPS k}{-}\widetilde\vv_k)
+(\widetilde\vv_k{\cdot}\nabla)(\vv_{\EPS k}{-}\widetilde\vv_k)\Big)
\\[-.1em]&\hspace*{8em}
-\frac12{\rm div}(\varrho_{\EPS k}\widetilde\vv_k)|\vv_{\EPS k}{-}\widetilde\vv_k|^2
\bigg)\,\d\xx\d t
+\int_\varOmega\frac{\varrho_0}2|\vv_0{-}\widetilde\vv_k(0)|^2\,\d \xx\ {\buildrel{k\to\infty}\over{\longrightarrow}}\ 0\,
\nonumber\end{align}
with some $\widetilde\nu_2>0$ depending on $p\ge2$, related to the inequality
$\widetilde\nu_2|G-\widetilde G|^p\le\nu_2(|G|^{p-2}G-|\widetilde G|^{p-2}\widetilde G)\Vdots(G-\widetilde G)$.
The term $\frac1k\vv_{\EPS k}{\cdot}(\vv_{\EPS k}{-}\widetilde\vv_k)$ converges
to 0 in $L^2(I{\times}\varOmega)$ because, due to \eq{est-rv2+},
$\|\frac1k\vv_{\EPS k}\|_{L^2(I\times\varOmega;\R^d)}^{}=\mathscr{O}(1/\sqrt{k})\to0$
while $\vv_{\EPS k}{-}\widetilde\vv_k$ is bounded in $L^2(I{\times}\varOmega;\R^d)$
thanks to \eq{basic-est-of-v}.
 Using \eq{v-conv} and that ${\rm div}(\varrho_{\EPS k}\widetilde\vv_k)=
\nabla\varrho_{\EPS k}{\cdot}\widetilde\vv_k+\varrho_{\EPS k}{\rm div}\widetilde\vv_k$
is bounded in $L^{\infty}(I;L^r(\varOmega))$ due to the already obtained bounds \eq{est-nabla-v} and \eq{est+rho}, we can see that
$\frac12{\rm div}(\varrho_{\EPS k}\widetilde\vv_k)|\vv_{\EPS k}{-}\widetilde\vv_k|^2\to0$ surely in $L^1(I{\times}\varOmega)$ since $r>2$. 

Thus (\ref{Euler-weak+}e,f) is proved.

The strong convergences \eq{Euler-weak+} allow for limit passage simply by
continuity in the Galerkin approximation of the regularized momentum equation
\eq{Euler1-weak-Galerkin}  to obtain  \eq{Euler1-thermo-finite+}, 
written without the regularizing term $\frac1k\vv_{\EPS k}$. Similarly, we
take the limit  in the regularized internal-energy equation \eq{Euler3-weak-Galerkin}
(or alternatively modified for \eq{Euler3-enthalpy-finite-reg})
 to obtain  \eq{Euler3-thermo-finite+},  written without the regularizing
term $\frac1k\pdt{}\theta_{\EPS k}$. Here, in addition to \cite{Roub24TVSE},
we use that the regularizing term $\frac1k\vv_{\EPS k}$ in \eq{Euler1-thermo-finite+}
converges to zero strongly in $L^2(I{\times}\varOmega;\R^d)$ due to the
latter estimate in \eq{est-rv2+} and also the mentioned convergence
$\frac1k\theta_{\EPS k}\to0$ used for the weak formulation \eq{Euler3-weak-Galerkin}.

\medskip\noindent{\it Step 5: Further a-priori estimates uniform for $\EPS>0$}.
First, let us emphasize that the estimates \eq{Euler-est}--\eq{est-nabla-v} and
\eq{est+}--\eq{est-adibatic-heat} are uniform not only with respect to $k\in\N$ but
also in $\EPS>0$. Thus they are inherited also for $(\vv_\EPS,\FF_\EPS,\theta_\EPS)$.
In particular, since we have $\theta_\EPS\ge0$ granted simply by construction of
the limit \eq{Euler-disc-sln-theta}, the estimate \eq{est-phi+} inherited for
$(\FF_\EPS,\theta_\EPS)$ with the assumption \eq{Euler-ass-psi-1} leads to
\begin{align}\label{Loo-est-for-theta}
\|\theta_\EPS\|_{L^\infty(I;L^{1+\ALPHA}(\varOmega))}^{}\le C\,.
\end{align}

We may benefit from  the fact  that the regularized internal-energy equation
is continuous. Therefore, one can employ the $L^1$-theory for heat-transfer-type
equations which is based on the ``nonlinear'' test by $\chi_\ZETA(\theta_\EPS)$
applied to the internal-energy 
equation \eq{Euler3-thermo-finite+} written for $(\vv_\EPS,\FF_\EPS,\theta_\EPS)$
instead of $(\vv,\FF,\theta)$ without the term $\pdt{}\theta/k$. Specifically,
we will use the nonlinear increasing function $\chi_\ZETA:[0,+\infty)\to[0,1]$
defined as
\begin{align}\label{test-chi}
\chi_\ZETA(\theta):=1-\frac1{(1{+}\theta)^\ZETA}\,,\ \ \ \ \ZETA>\ALPHA\,,
\end{align}
as also used for $\ALPHA=0$ in \cite{Naum06EWSE},
simplifying the original idea of L.\,Boccardo, T.\,Gallou\"et, et al.\ 
\cite{BDGO97NPDE,BocGal89NEPE} and extending it for the mechanically coupled
systems. Cf.\ also \cite{BuFeMaNSFS,FeiMal06NSET} where
other functions were used, namely $\chi_\ZETA(\theta):=\theta^\ZETA$
or $\chi_\ZETA(\theta):=1/\theta^\ZETA$ with $0<\ZETA<1$.
Importantly, here we have $\chi_\ZETA(\theta_\EPS(t,\cdot))\in H^1(\varOmega)$,
so it is a  admissible  test function, because 
$0\le\theta_\EPS(t,\cdot)\in H^1(\varOmega)$ has already been proved
and because $\chi_\ZETA$ is Lipschitz continuous on $[0,+\infty)$. 

The rather technical arguments in \cite{Roub24TVSE} for $\ALPHA=0$ 
have to be generalized here for $\ALPHA\ge0$, also combining the interpolation
from \cite[Sect.8.2]{KruRou19MMCM}. Since this generalization is not
entirely straightforward at some spots, we will present this part quite in
detail. We consider $1\le \EXP<2$ and estimate
the $L^\EXP$-norm  of $\nabla\theta_\EPS$ by H\"older's inequality as 
\begin{align}\label{8-**-+}
&\!\!\int_0^T\!\!\!\int_\varOmega\!|\nabla\theta_\EPS|^\EXP\,\d\xx\d t
\le C_1\bigg(\underbrace{\int_0^T\!\!\!
\big\|1{+}\theta_\EPS(t,\cdot)\big\|^{(1+\ZETA)\EXP/(\TWO-\EXP)}
_{L^{(1+\ZETA)\EXP/(\TWO-\EXP)}(\varOmega)}\,\d t}_{\displaystyle\ \ \ =:I_{\EXP,\ZETA}^{(1)}(\theta_\EPS)}\bigg)^{1-\EXP/\TWO}
\bigg(\underbrace{\int_0^T\!\!\!\int_\varOmega\!
\chi_\ZETA'(\theta_\EPS)|\nabla\theta_\EPS|^\TWO}
_{\displaystyle\ \ \ =:I_{\ZETA}^{(2)}(\theta_\EPS)}\bigg)^{\EXP/\TWO}\hspace*{-.5em}
\end{align}
with $\chi_\ZETA$ from \eq{test-chi} so that
$\chi_\ZETA'(\theta)=\ZETA/(1{+}\theta)^{1+\ZETA}$ and with a constant $C_1$ dependent
on $\ZETA$ and $\EXP$. Then we interpolate the Lebesgue space
$L^{(1+\ZETA)\EXP/(\TWO-\EXP)}(\varOmega)$ between  $W^{1,\EXP}(\varOmega)$ and 
$L^{1+\ALPHA}(\varOmega)$ to exploit \eq{Loo-est-for-theta}.
More specifically, by the Gagliardo-Nirenberg inequality, we obtain
\begin{align}\label{8-cond}
\big\|1{+}\theta_\EPS(t,\cdot)\big\|_{L^{(1+\ZETA)\EXP/(\TWO-\EXP)}(\varOmega)}^{}
\le C_2&\Big(1+\big\|\nabla\theta_\EPS(t,\cdot)\big\|_{L^\EXP(\varOmega;\R^d)}\Big)^\lambda
\\
&\text{ for }\ \ \frac{2-\EXP}{(1{+}\ZETA)\EXP}\ge\lambda\Big(\frac1\EXP-\frac1d\Big)+\frac{1{-}\lambda}{1{+}\ALPHA}\ \text{ with }\ 0<\lambda\le1\,.
\nonumber\end{align}
Thus, we obtain
$I_{\EXP,\ZETA}^{(1)}(\theta_\EPS)\le C_3(1+\int_0^T\!\int_\varOmega\big|\nabla\theta_\EPS\big|^\EXP\,\d\xx\d t)$
providing
\begin{align}\label{choice-lambda}
\lambda=\frac{\TWO{-}\EXP}{1{+}\ZETA}\,,
\end{align}
cf.\ also \cite[Formulas (8.2.14)-(8.2.16)]{KruRou19MMCM}.
Combining it with \eq{8-**-+}, we obtain
\begin{align}\|\nabla\theta_\EPS\|_{L^\EXP(I\times\varOmega;\R^d)}^\EXP=C_1C_3\big(1+\|\nabla\theta_\EPS\|_{L^\EXP(I\times\varOmega)}^\EXP\big)^{1-\EXP/2}I_{\ZETA}^{(2)}(\theta_\EPS)^{\EXP/2_{_{_{}}}}_{^{^{^{}}}}\,.
\label{8-***}
\end{align}
Furthermore, we are to estimate $I_{\ZETA}^{(2)}(\theta_\EPS)$ in \eq{8-**-+}.
Like \eq{primitive-theta-c}, we now use a primitive function ${\mathcal X}_\ZETA$ to
$\theta\mapsto\chi_\ZETA(\theta)\ENG_\theta'(\Fe,\theta)$ depending
smoothly on $\Fe$, specifically
\begin{align}
{\mathcal X}_\ZETA(\Fe,\theta)
=\int_0^1\!\!\theta\chi_\ZETA(r\theta)\ENG_\theta'(\Fe,r\theta)\,\d r\,.
\label{primitive+}\end{align}
To modify \eq{heat-test-theta} appropriately, we use now the calculus
\begin{align}\label{Euler-thermodynam3-test--}
\!\!\!\int_\varOmega\!\chi_\ZETA(\theta)\pdt\Eng\,\d\xx
&=\!\int_\varOmega\!\chi_\ZETA(\theta)\ENG_\theta'(\Fe,\theta)\pdt{\theta}
+\chi_\ZETA(\theta)\ENG_\Fe'(\Fe,\theta){:}\pdt{\Fe\!}\,\,\d\xx
\\&\nonumber=\frac{\d}{\d t}\int_\varOmega{\mathcal X}_\ZETA(\Fe,\theta)\,\d\xx
-\int_\varOmega\big[{\mathscr X}_\ZETA\big]_\Fe'(\FF,\theta){:}\pdt{\Fe}\,\d\xx
\\&\qquad\qquad\qquad\text{ where }\
{\mathscr X}_\ZETA(\FF,\theta):={\mathcal X}_\ZETA(\Fe,\theta)
-\chi_\ZETA(\theta)\ENG(\Fe,\theta)\,.
\nonumber\end{align}
In view of \eq{primitive+}, it holds $[{\mathscr X}_\ZETA]_\Fe'(\Fe,\theta)
=\int_0^1\theta\chi_\ZETA(r\theta)\ENG_{\Fe\theta}''(\Fe,r\theta)\,\d r
-\chi_\ZETA(\theta)\ENG_{\Fe}'(\Fe,\theta)$. Altogether, testing
\eq{Euler3-thermo-finite+} with \eq{Euler-thermodynam-BC-2-reg} by
$\chi_\ZETA(\theta_\EPS)$ gives
\begin{align}\label{Euler-thermodynam3-test+++}
\frac{\d}{\d t}\int_\varOmega&\!{\mathcal X}_\ZETA(\FFeps,\theta_\EPS)\,\d\xx
+\!\int_\varOmega
\chi_\ZETA'(\theta_\EPS)\kappa(\FFeps,\theta_\EPS)|\nabla\theta_\EPS|^2\,\d\xx
\\&\nonumber=\!\int_\varOmega\!\bigg(
\frac{\nu_1|\ee(\vv_\EPS)|^2+\nu_2|\nabla^2\vv_\EPS|^p}
{1{+}\EPS|\ee(\vv_\EPS)|^q{+}\EPS|\nabla\ee(\vv_\EPS)|^p\!}
\,\chi_\ZETA(\theta_\EPS)
-\ENG(\FFeps,\theta_\EPS)\chi_\ZETA'(\theta_\EPS)\vv_\EPS{\cdot}\nabla\theta_\EPS
\\&\qquad
+\big[{\mathscr X}_\ZETA\big]_\Fe'(\FFeps,\theta_\EPS){:}\pdt{\FFeps\!\!}
+\chi_\ZETA(\theta_\EPS)\TT(\FFeps,\theta_\EPS){:}\ee(\vv_\EPS)
\!\bigg)\d\xx+\!\int_\varGamma\!h_\EPS\chi_\ZETA(\theta_\EPS)\d S.
\nonumber\end{align}
We realize that $\chi_\ZETA'(\theta)=\ZETA/(1{+}\theta)^{1+\ZETA}$ as used
already in \eq{8-**-+} and that ${\mathcal X}_\ZETA(\FFeps,\theta_\EPS)\ge
c_K\theta_\EPS^{1+\ALPHA}$ with some $c_K$ 
due to \eq{Euler-ass-psi-1+++}; again $K$ is a compact subset of ${\rm GL}^+(d)$
related here with the already proved estimates \eq{est+Fes} inherited for $\FF_\EPS$.
The convective term in \eq{Euler-thermodynam3-test+++} can be estimated,
for any $\delta>0$, as
\begin{align}\label{est-of-convective}
\int_\varOmega\Eng_\EPS\chi_\ZETA'(\theta_\EPS)\vv_\EPS{\cdot}\nabla\theta_\EPS\,\d\xx
&\le\frac1\delta\int_\varOmega\chi_\ZETA'(\theta_\EPS)|\vv_\EPS|^2\Eng_\EPS^2\,\d\xx
+\delta\int_\varOmega\chi_\ZETA'(\theta_\EPS)|\nabla\theta_\EPS|^\TWO\,\d\xx
\\&=\frac1\delta\int_\varOmega\chi_\ZETA'(\theta_\EPS)|\vv_\EPS|^2
\Eng_\EPS^2\,\d\xx+\delta I_{\ZETA}^{(2)}(\theta_\EPS)\,.
\nonumber\end{align}
For $c_K$ acting in \eq{Euler-ass-kappa}, using \eq{Euler-thermodynam3-test+++}
integrated over $I=[0,T]$, we further estimate:
\begin{align}\label{+++}
I_{\ZETA}^{(2)}(\theta_\EPS)
&\le\frac1{\ZETA c_K}\int_0^T\!\!\!\int_\varOmega\!\kappa(\FFeps,\theta_\EPS)
\nabla\theta_\EPS{\cdot}\nabla\chi_\ZETA(\theta_\EPS)\,\d\xx
\\&\nonumber
\le\frac1{\ZETA c_K}\int_0^T\!\!\!\int_\varOmega\!\kappa(\FFeps,\theta_\EPS)
\nabla\theta_\EPS{\cdot}\nabla\chi_\ZETA(\theta_\EPS)\,\d\xx
+\frac1{\ZETA c_K}\int_\varOmega\!{\mathcal X}_\ZETA(\FFeps(T),\theta_\EPS(T))\,\d\xx
\\&\nonumber
\!\!\stackrel{\scriptsize\eq{Euler-thermodynam3-test+++}}{=}\!\!
\frac1{\ZETA c_K}\bigg(\int_\varOmega\!{\mathcal X}_\ZETA(\Fe_0,\theta_{0,\EPS})\,\d\xx
+\!\int_0^T\!\!\!\int_\varOmega\!
\bigg(\frac{\nu_1|\ee(\vv_\EPS)|^2
+\nu_2|\nabla^2\vv_\EPS|^p}{1{+}\EPS|\ee(\vv_\EPS)|^q{+}\EPS|\nabla\ee(\vv_\EPS)|^p\!}
\,\chi_\ZETA(\theta_\EPS)
\\&\nonumber\qquad
-\ENG(\FFeps,\theta_\EPS)\chi_\ZETA'(\theta_\EPS)\vv_\EPS{\cdot}\nabla\theta_\EPS
+\big[{\mathscr X}_\ZETA\big]_\Fe'(\FFeps,\theta_\EPS){:}\pdt{\FFeps}
\\[-.4em]&\nonumber\qquad\qquad\qquad\qquad
+
\TT(\FFeps,\theta_\EPS){:}\ee(\vv_\EPS)\chi_\ZETA(\theta_\EPS)
\bigg)\,\d\xx\d t
+\!\int_0^T\!\!\!\int_\varGamma\!h_\EPS\chi_\ZETA(\theta_\EPS)\,\d S\d t\bigg)
\\[-.2em]&\nonumber
\le
\frac1{\ZETA c_K}\bigg(
\big\|{\mathcal X}_\ZETA(\Fe_0,\theta_{0,\EPS})\big\|_{L^1(\varOmega)}\!
+\big\|\nu_1|\ee(\vv_\EPS)|^2+\nu_2|\nabla^2\vv_\EPS|^p\big\|_{L^1(I\times\varOmega)}\!
+\big\|h\big\|_{L^1(I\times\varGamma)}\!
\\&\quad\nonumber
+\big\|\TT(\FFeps,\theta_\EPS){:}\ee(\vv_\EPS)\big\|_{L^1(I\times\varOmega)}\!
+\!\!\int_0^T\!\!\!\big\|\big[{\mathscr X}_\ZETA\big]_\Fe'(\FFeps,\theta_\EPS)
\big\|_{L^{r'}(\varOmega;\R^{d\times d})}^{r'}\!
+\Big\|\pdt{\FFeps}\Big\|_{L^r(\varOmega;\R^{d\times d})}^r\d t\!\!\!\!\!
\\[-.2em]&\quad
+\frac1\delta\|\vv_\EPS\|_{L^2(I;L^\infty(\varOmega;\R^d))}^{\TWOprime}
\|\chi_\ZETA'(\theta_\EPS)\ENG^{\TWOprime}(\FFeps,\theta_\EPS)
\|_{L^\infty(I;L^1(\varOmega))}^{}\!\bigg)
+\frac\delta{c_K}I_{\ZETA}^{(2)}(\theta_\EPS)\,.
\nonumber\end{align}
Here we used the bound of $\pdt{}\FF_\EPS$ in
$L^1(I;L^r(\varOmega;\R^{d\times d}))$, as follows from due to
\eq{Euler-disc-sln-F} with the exponent $r$ arbitrarily large  because
of the qualification $\FF_0\in W^{1,r}(\varOmega;\R^{d\times d})$ for any $r$,
cf.\ \eq{Euler-ass-Fe0}.  By the qualification \eq{Euler-ass-psi-1+++},
we have $|[{\mathscr X}_\ZETA]_\Fe'(\FF,\theta)|\le C(1{+}\theta^{1+\ALPHA})$,
so that we need the estimation of the term
\begin{align}\label{++++}
\|[{\mathscr X}_\ZETA]_\Fe'(\FFeps,\theta_\EPS)\|_{L^{r'}(\varOmega;\R^{d\times d})}^{r'}
&\le C\|1{+}\theta_\EPS\|_{L^{r'(1+\ALPHA)}(\varOmega)}^{r'(1+\ALPHA)}
\\&\le
C_{r,\ALPHA}\big(1+\|\nabla\theta_\EPS\|_{L^\EXP(\varOmega;\R^d)}^b\,\big)
\le C_{r,\ALPHA,\delta}+\delta\|\nabla\theta_\EPS\|_{L^\EXP(\varOmega;\R^d)}^\EXP
\nonumber\end{align}
with an arbitrarily small $b>0$ and thus also arbitrary small $\delta>0$.
Here we used  the already obtained estimate \eq{Loo-est-for-theta}
and the Gagliardo-Nirenberg inequality to interpolate
$L^{r'(1+\ALPHA)}(\varOmega)$ between $L^{1+\ALPHA}(\varOmega)$ and
$W^{1,\EXP}(\varOmega)$  for any $r'>1$. We can take  $r'>1$ arbitrarily
small  as $r$ cn be taken arbitrarily large as mentioned above..
Due to \eq{growth-T}, we can estimate the power of the conservative
stress $\TT(\FFeps,\theta_\EPS){:}\ee(\vv_\EPS)$  from  \eq{+++} in
$L^\infty(I;L^1(\varOmega))$. The penultimate term in \eq{+++} is a-priori bounded
independently of $\EPS$ for $\ZETA>\ALPHA$ fixed because
$\ENG(\FF,\theta)=\mathscr{O}(\theta^{1+\ALPHA})$ due to \eq{Euler-ass-psi-1+++}
and  because  $\chi_\ZETA'(\theta)=\mathscr{O}(1/\theta^{1+\ALPHA})$ uniformly for
$\ZETA>\ALPHA$. Let us emphasize that here we rely on  the fact  that $\ZETA>\ALPHA$
as stated  in \eq{test-chi}. As a result, we have
$[\chi_\ZETA'(\cdot)\ENG^\TWO(\FF,\cdot)](\theta)=\mathscr{O}(\theta^{1+\ALPHA})$.
Thus the estimate \eq{Loo-est-for-theta} guarantees
$\chi_\ZETA'(\theta_\EPS)\ENG^{\TWOprime}(\FFeps,\theta_\EPS)$ bounded
in $L^\infty(I;L^1(\varOmega))$ while $|\vv_\EPS|^{\TWOprime}$
is surely bounded in $L^1(I;L^\infty(\varOmega))$, cf.\ 
 the estimate  \eq{basic-est-of-v} inherited for $\vv_\EPS$ independently
of $\EPS>0$. Eventually, when choosing $\delta<c_K$, we can absorb the last
term in the left-hand side. Similarly, combining it with \eq{8-***}, also the
 last term in \eq{++++} integrated over the time interval $I$ 
can be absorbed
in the left-hand side of \eq{8-***}.

Substituting $\lambda
$ from \eq{choice-lambda} into \eq{8-cond}, by some algebra
we obtain the bound $\EXP\le(d{+}2{+}2\ALPHA{-}\ZETA d)/(d{+}1{+}\ALPHA)$,
cf.\ also \cite[Formula (8.2.18)]{KruRou19MMCM}. Taking into account
$\ZETA>\ALPHA$ in \eq{test-chi}, we obtain the restriction on the integrability
exponent for the temperature gradient:
\begin{align}
\EXP<\frac{d+2+(2{-}d)\ALPHA}{d+1+\ALPHA}\,.
\label{bound-for-mu}\end{align}
Simultaneously, $\EXP\ge1$ is desirable, which needs $\ALPHA<1/2$ if $d=3$
or $\ALPHA<1$ if $d=2$ as assumed in (\ref{Euler-ass}e,f). Realizing the embedding
$L^\infty(I;L^{1+\ALPHA}(\varOmega))\,\cap\,L^\EXP(I;W^{1,\EXP}(\varOmega))\subset
L^s(I{\times}\varOmega)$ with $1\le s<\EXP(1+(1{+}\ALPHA)/d)=(d+2+(2{-}d)\ALPHA)/d$,
in total we have proved
\begin{subequations}\label{est-W-eps}\begin{align}
&\label{est-theta-1}
\big\|\theta_{\EPS}\big\|_{L^\infty(I;L^{1+\ALPHA}(\varOmega))\,\cap\,L^{ s}(I{\times}\varOmega)}^{}\le C
\ \ \ \text{ with }\ 1\le s< \smaller{\frac{d{+}2{+}(2{-}d)\ALPHA}d}
\ \ \ \text{ and }\ \ \
\\[-.3em]&\label{est-theta-2}
\big\|\nabla\theta_{\EPS}\big\|_{L^\EXP(I\times\varOmega;\R^d)}^{}\le C\ \ \
\text{ with }\ 1\le \EXP<\smaller{\frac{d{+}2{+}(2{-}d)\ALPHA}{d{+}1{+}\ALPHA}}\,.
\intertext{Exploiting the calculus
$\nabla\Eng_{\EPS}=\ENG_\theta'(\FF_{\EPS},\theta_{\EPS})\nabla\theta_{\EPS}+
\ENG_\FF'(\FF_{\EPS},\theta_{\EPS})\nabla\FF_{\EPS}$
with $\nabla\FF_{\EPS}$ bounded in $L^\infty(I;L^r(\varOmega;\R^{d\times d\times d}))$
for and $1\le r<+\infty$
and relying on the assumption \eq{Euler-ass-psi-1+++}, we have also the
bound on $\nabla\Eng_{\EPS}$ in $L^\EXP(I;L
^{\EXP(1+\ALPHA)/(1+\ALPHA+\ALPHA\EXP)}(\varOmega;\R^d))$, so
that}
&\|\Eng_\EPS\|_{L^\infty(I;L^1(\varOmega))\,\cap\,L^\EXP(I;W^{1,\EXP(1+\ALPHA)/(1+\ALPHA+\ALPHA\EXP)}(\varOmega))}^{}\le C_\EXP\ \ \text{ with $\ \mu\ $ from \eq{est-theta-2}}\,.
\intertext{Since $\Eng_\EPS=\ENG(\FF_\EPS,\theta_\EPS)$ and $\ENG(\FF,\cdot)$
is bounded from below due to \eq{ass-stress-control} and has at most
$(1{+}\ALPHA)$-polynomial growth \eq{growth-int-eng}, uniformly for $\FF\in K$,
we also have}
\label{est-theta-4}
&\|\Eng_\EPS\|_{L^{s/(1+\ALPHA)}I{\times}\varOmega)}^{}\le C\ \ \text{ with $\ s\ $ from \eq{est-theta-1}}\,.
\end{align}\end{subequations}

\medskip\noindent{\it Step 6: Convergence for $\EPS\to0$}. 
Using the Banach selection principle as in Step~4, now also taking
the estimates \eq{Euler-est}--\eq{est-nabla-v} and
\eq{est+}--\eq{est-adibatic-heat} inherited for $(\vv_\EPS,\FF_\EPS,\theta_\EPS)$
and \eq{est-W-eps} into account, we can extract some subsequence of
$\{(\varrho_{\EPS},\vv_{\EPS},\FFeps,\theta_{\EPS},\Eng_{\EPS})\}_{\EPS>0}$
and its limit
$(\varrho,\vv,\FF,\theta,\Eng):I\to H^1(\varOmega)\times
L^2(\varOmega;\R^d)\times H^1(\varOmega;\R^{d\times d})\times
L^{1+\ALPHA}(\varOmega)\times L^1(\varOmega)$
such that, for any $1\le r<+\infty$,
\begin{subequations}\label{Euler-weak++}
\begin{align}
&\!\!\varrho_{\EPS}\to\varrho\!\!\!&&\text{weakly* in $\
L^\infty(I;W^{1,r}(\varOmega))\,\cap\,W^{1,
p}(I;L^r(\varOmega))$}
\\[-.3em]&&&\hspace*{12em}\text{and strongly in $C(I{\times}\barOmega)$\,,}
\nonumber\\\label{Euler-weak-sln-v}
&\!\!\vv_{\EPS}\to\vv\!\!\!&&\text{weakly* in $\ L^\infty(I;L^2(\varOmega;\R^d))\cap
L^{
p}(I;W^{2,p}(\varOmega;\R^d))$\,,}&&
\\\label{Euler-weak-F}
&\!\!\FFeps\to\FF
\!\!\!\!\!\!&&\text{weakly* in $\ L^\infty(I;W^{1,r}(\varOmega;\R^{d\times d}))\,\cap\,
W^{1,
p}(I;L^r(\varOmega;\R^{d\times d}))$}\!\!
\\[-.2em]&&&\hspace*{12em}\text{and strongly in
$C(I{\times}\barOmega;\R^{d\times d})$\,,}
\nonumber\\\label{Euler-weak-sln-theta}
&\!\!\Eng_{\EPS}\to\Eng\!\!\!\!\!\!&&\text{weakly in $\
L^\EXP(I;W^{1,\EXP(1+\ALPHA)/(1+\ALPHA+\ALPHA\EXP)}(\varOmega))$\,, and}\!\!
\\
&\!\!\theta_{\EPS}\to\theta\!\!\!&&
\text{weakly in $\ L^\EXP(I;W^{1,\EXP}(\varOmega))\ $ with $\ \mu\ $
from \eq{est-theta-2}.}
\label{Euler-weak-sln-w++}
\end{align}\end{subequations}
We now also have the estimate of $\pdt{}\Eng_{\EPS}$ in $L^1(I;H^3(\varOmega)^*)$
relying on $H^3(\varOmega)\subset W^{1,\infty}(\varOmega)$ for $d\le3$. So, 
like \eq{w-conv} and using \eq{est-theta-4}, by the Aubin-Lions theorem we now
have
\begin{subequations}\label{Euler-weak+++}
\begin{align}
&\Eng_\EPS\to\Eng\!\!\!&&\hspace*{-5em}\text{strongly in $L^{s/(1{+}\ALPHA)}(I{\times}\varOmega)$\ \ with $\ s\ $ from \eq{est-theta-1}}\,.
\intertext{Then, we realize that $\Ent_\EPS=\ENT(\FF_\EPS,\theta_\EPS)$ 
with $\Ent_\EPS=\Eng_\EPS-\ENG(\FF_\EPS,0)$ converges to
$\Ent=\Eng-\ENG(\FF,0)=\ENT(\FF,\theta)$ strongly in $L^{s/(1{+}\ALPHA)}(I{\times}\varOmega)$ and, like in \eq{z-conv}, we use again continuity of
$(\FF,\Ent)\mapsto[\ENT_{\text{\sc m}}^{}(\FF,\cdot)]^{-1}(\Ent)$ and realize that
actually
$[\ENT_{\text{\sc m}}^{}(\FF,\cdot)]^{-1}(\Ent_\EPS)=[\ENT(\FF,\cdot)]^{-1}(\Ent_\EPS)$
since $\Ent_\EPS\ge0$. By this arguments, we also have}
&\theta_\EPS\to\theta=[\ENT(\FF,\cdot)]^{-1}(\Ent)\!\!\!&&
\hspace*{.1em}\text{strongly in $\ L^s(I{\times}\varOmega)
\ $ with $\ s\ $ again from \eq{est-theta-1}.}
\intertext{By the continuity of $\TT(\cdot,\cdot)$ and of $\kappa$, we have also}
&\TT(\FFeps,\theta_\EPS)\to\TT(\FF,\theta)
\hspace*{.5em}&&\hspace*{-1.3em}\text{strongly in
$L^{r}(I{\times}\varOmega;\R_{\rm sym}^{d\times d})$ for any $1\le r<
\dfrac{d{+}2{+}(2{-}d)\ALPHA\!}{d(1{+}\ALPHA)}
$\,, and}
\label{Euler-weak-stress}
\\[-.3em]&\label{m-strongly+}
\kappa(\FFeps,\theta_\EPS)\to\kappa(\FF,\theta)
&&\hspace*{-1.3em}\text{strongly in $L^r(I{\times}\varOmega)$ for any $1\le r<\infty$.}
\end{align}\end{subequations}
For the exponent in \eq{Euler-weak-stress}, we used the estimate \eq{est-theta-1}
together with the energy-controlled-stress condition \eq{ass-stress-control}
exploiting also \eq{growth-int-eng}. Also the strong convergence of velocity
gradients like  (\ref{Euler-weak+}e,f) can be proved for $\vv_\EPS$ instead of
$\vv_{\EPS k}$ analogously as in Step~4.

The limit passage in the weak formulation of the system
\eq{Euler-thermo-finite+} without the terms $\frac1k\vv$ and
$\frac1k\pdt{}\theta$ towards the solution of the system
\eq{Euler-thermo-finite} according Definition~\ref{def} is then easy by
continuity and convergence $h_\EPS\to h$ in $L^1(I{\times}\varGamma)$.

\medskip\noindent{\it  Step 7: Energy balances}.
In fact, any solution according Definition~\ref{def} that possesses
the regularity as specified in Theorem~\ref{prop-Euler} satisfies
the energy balances \eq{Euler-engr-finite-}, \eq{Euler-engr-finite},
and \eq{Euler-engr-finite--} integrated over a current time interval $[0,t]$.
It is important that the above tests of the four equations in
\eq{Euler-thermo-finite} successively by $\frac12|\vv|^2$, $\vv$, $\phi'(\FF)$,
and 1 are analytically legitimate; here we again refer to \cite{Roub24TVSE}.
\end{proof}

\begin{remark}[{\sl Fourier boundary conditions}]\label{rem-Fourier}\upshape
The last boundary condition \eq{Euler-thermodynam-BC}
can be generalized towards the Robin-type (in the heat context also called
Fourier) condition $\kappa(\FF,\theta)\nabla\theta{\cdot}\nn+b\theta=h$
with $b\in L^\infty(\varGamma)$ non-negative. Then \eq{Euler-thermodynam-BC-2-reg}
should be modified as $\kappa(\FF,\theta)\nabla\theta{\cdot}\nn+b\theta^+=h_\EPS$,
which would  again  work
 also  in the Galerkin approximation where also the negative
values of temperature may occur. The limit passage in Step~6 would
then rely on the compactness of the trace operator
$W^{1,\EXP}(\varOmega)\to L^1(\varGamma)$ with $\EXP>1$ from \eq{bound-for-mu}.
\end{remark}

\begin{remark}[{\sl  Non-negativity of the temperature}]\upshape
 Proving $\theta\ge0$  requires the test 
of the internal-energy equation \eq{Euler3-thermo-finite} by $\theta^-$
which is not in duality with $\DT\Eng$ in general, however. Thus 
we should realize that, rigorously, we have proved this non-negativity only
for at least one solution  constructed by the semi-Galerkin procedure
above.  Finally, note that the entropy balance \eq{Euler-entropy-finite}
needs positivity of $\theta$, which is not obvious and likely very nontrivial
to prove even for the non-negative temperature mentioned above, because 
the heat capacity in our Eulerian formulation inevitably also depends on $\FF$,
which complicates the usual arguments.
\end{remark}

\begin{remark}[{\sl Smoothness of $\psi$}]\label{rem-smooth}\upshape
Noteworthy, the assumption \eq{Euler-ass-psi-1}   has needed 
existence of $\psi_{\theta\theta}''$ and \eq{Euler-ass-psi-1+++}  has
needed   $\psi_{\FF\theta}''$ and even $\psi_{\FF\theta\theta}'''$,  while the system \eq{Euler-thermo-finite} itself
involves only $\psi_{\FF}'$ and $\psi_{\theta}'$. In fact,  a  more careful
formulation of the assumptions (\ref{Euler-ass}e,f) and a suitable mollification
of such less smooth $\psi$ would allow  the use of 
Theorem~\ref{prop-Euler} followed by another limit passage  to admit 
a discontinuous $\psi_{\theta\theta}''$.  This could make it
possible  to describe 2nd-order phase transitions  of the
Stefan type, during which $\ENT(\FF,\cdot)$ may have a jump
at a transition temperature. 
\end{remark}

\begin{remark}[{\sl Exploiting the entropy balance -- an open problem}]
\label{rem-entrop}\upshape
One may be tempted to use the entropy equation \eq{Euler-entropy-finite}
for a direct estimation of the temperature gradient in order to avoid,
like \cite{FeiNov09SLTV} does, the use of the technically demanding
$L^1$-theory for the heat equation based on \cite{BDGO97NPDE,BocGal89NEPE}.
We recall that \eq{Euler-entropy-finite} was obtained by using the test
of \eq{Euler3-thermo-finite} by $1/\theta$, i.e.\ by a so-called {\it coldness}.
However, such nonlinear test is incompatible with a Galerkin approximation,
so one is tempted to formulate the system directly in terms of the coldness,
following the idea of Truesdell \cite{True69RT} who articulated that `in
continuum mechanics the reciprocal of the temperature, which may be called the
coldness $\vartheta$, is often more convenient for the mathematical theory'.
This would potentially allow for omitting the intermediate test by $\theta$
in Step~3 and the technically demanding $L^1$-theory for estimating
$\nabla\theta$ in Step~5 in the above proof. In terms of $\vartheta:=1/\theta$,
one can define $\wh\psi(\FF,\vartheta):=\psi(\FF,1/\vartheta)$, 
$\wh\ENG(\FF,\vartheta):=\ENG(\FF,1/\vartheta)$, 
$\wh\eta(\FF,\vartheta):=\eta(\FF,1/\vartheta)$, 
$\wh\kappa(\FF,\vartheta):=\kappa(\FF,1/\vartheta)$,
and the Cauchy stress $\wh\TT(\FF,\vartheta):=\wh\psi_\FF'(\FF,\vartheta)\FF^\top\!\!
+\wh\psi(\FF,\vartheta)\bbI$. The relation between $\wh\ENG$ and $\wh\psi$ is
\begin{align}\nonumber
\wh\ENG(\FF,\vartheta):=\ENG\Big(\FF,\frac1\vartheta\Big)&=
\psi\Big(\FF,\frac1\vartheta\Big)
-\frac1\vartheta\psi_\theta'\Big(\FF,\frac1\vartheta\Big)=
\wh\psi(\FF,\vartheta)
-\frac1\vartheta\psi_\theta'\Big(\FF,\frac1\vartheta\Big)
\\\nonumber&=
\wh\psi(\FF,\vartheta)
-\vartheta\frac1{\vartheta^2}\psi_\theta'\Big(\FF,\frac1\vartheta\Big)
=
\wh\psi(\FF,\vartheta)
+\vartheta\wh\psi_\vartheta'(\FF,\vartheta)\,.
\end{align}
The heat equation in terms of the internal energy \eq{Euler3-thermo-finite} is
then transformed to
\begin{align}\nonumber&
\DT\Eng={\rm div}\big(\wh\kappa(\FF,\vartheta)\nabla\vartheta\big)-({\rm div}\,\vv)\Eng
-\nu_1|\ee(\vv)|^2-\nu_2|\Nabla^2\vv|^p-\wh\TT(\FF,\vartheta)\Colon\Nabla\vv
\end{align}
with $\Eng=\wh\ENG(\FF,\vartheta)$.
This equation well respects physics and can be discretized by the Galerkin method,
allowing the test by $1$ to obtain the total-energy equation (from which one can
obtain $\vartheta\ge0$ even in the discrete approximation
due to the blowup of $\wh\ENG(\FF,\cdot)$ at $\vartheta=0$ and also
an estimate for $\wh\eta\,$) and by $\vartheta$ relying on non-negativity of
entropy, from which one can get an $L^2$-estimate of
$\sqrt{\wh\kappa(\FF,\vartheta)}\nabla\vartheta$ and, in some sense, also the
dissipation rate. Then, from the mechanical-energy balance as
\eq{Euler-engr-finite-}, one can improve the estimate of the dissipation rate.
This seems a promissing strategy but, having in mind Example~\ref{exa-neo} below,
$\wh\ENG(\FF,\cdot)$ has a singularity like $\vartheta^{-1-\ALPHA}$ so that estimation
of $\nabla\Eng=\wh\ENG_\FF'(\FF,\vartheta){:}\nabla\FF
+\wh\ENG_\vartheta'(\FF,\vartheta)\nabla\vartheta$ calls for a singularity of
$\wh\kappa(\FF,\cdot)$ like $\vartheta^{-3-\ALPHA}$, which eventually causes problems
with convergence of the Galerkin approximation.
\end{remark}

\section{Examples and remarks}\label{sec-exa}

\def\Eref{{{\rm E}}}

Let us complete this article with physically relevant examples that comply
with the assumptions \eq{Euler-ass}. We will formulate them in terms of the
referential free energy $\uppsi$ as used in Remark~\ref{rem-engineering},
which will also be defined for negative temperature for the sake of the
Galerkin method as used in the proof of Theorem~\ref{prop-Euler}. 
In general, such that $\uppsi(\FF,\cdot)$ is to be concave not to
conflict with non-negativity of the heat capacity.
We distinguish the actual free energy, the actual  internal energy, and
the actual heat capacity, i.e.\ respectively
\begin{align}\label{referential-psi}
\psi(\FF,\theta)=\frac{\uppsi(\FF,\theta)}{\det\FF}\,,\ \ \ \ \
\ENG(\FF,\theta)=\frac{\Eref(\FF,\theta)}{\det\FF}\,,\ \ \text{ and }\ \
c(\FF,\theta)=\frac{{\rm c}(\FF,\theta)}{\det\FF}
\end{align}
for $\det\FF>0$. In \eq{referential-psi}, we have used the referential internal
energy $\Eref(\FF,\theta)=\uppsi(\FF,\theta)+\theta\upeta(\FF,\theta)$ with
the referential entropy $\upeta(\FF,\theta)=-\uppsi_\theta'(\FF,\theta)$,
and the referential heat capacity
${\rm c}(\FF,\theta)=\theta\upeta_\theta'(\FF,\theta)
=-\theta\uppsi_{\theta\theta}''(\FF,\theta)$.

\begin{example}[{\sl A free energy of the neo-Hookean type}]\label{exa-neo}
\upshape
Denoting by $K_\text{\sc e}^{}$ and $G_\text{\sc e}^{}$ the elastic
bulk and shear moduli, respectively, $c_{\text{\sc v}}>0$ a referential
heat capacity at a reference temperature $\theta_0>0$ at constant volume, and
the exponent $\ALPHA>0$, a physically relevant example is
\begin{align}\label{free-neo-Hookean}
\uppsi(\FF,\theta)=\upphi(\FF)-\frac{c_{\text{\sc v}}\theta}{\!\ALPHA(1{+}\ALPHA)}\Big(\frac{\theta^+}{\theta_0}\Big)^\ALPHA
\ \ \text{ with }\ \ \upphi(\FF)=\frac12K_\text{\sc e}^{}\big(\det\FF{-}1\big)^2\!+
G_\text{\sc e}^{}\frac{{\rm tr}(\FF\FF^\top)}{(\det\FF)^{2/d}}
\,.
\end{align}
The actual entropy $\eta(\FF,\theta)=\upeta(\FF,\theta)/\!\det\FF=-\psi_\theta'(\FF,\theta)$ is then
\begin{align}\label{exa-ent}
\eta(\FF,\theta)=\frac{c_{\text{\sc v}}}{\ALPHA\det\FF}\Big(\frac{\theta^+}{\theta_0}\Big)^\ALPHA\,.
\end{align}
In particular, $\eta(\FF,0)=0$, i.e.\ it complies with the 3rd law of thermodynamics,
which states  that entropy is zero and, in particular, independent of the mechanical
state at zero temperature. Furthermore, the actual heat capacity
$c(\FF,\theta)=\theta\eta_\theta'(\FF,\theta)$ is
\begin{align}\label{exa-heat-cap}
c(\FF,\theta)=
\frac{c_{\text{\sc v}}}{\det\!\FF}\Big(\frac{\theta^+}{\theta_0}\Big)^\ALPHA
\end{align}
and obviously vanishes at zero temperature.
Note that we defined $\uppsi$ also for negative temperatures (cf.\
Figure~\ref{fig-free-energy}-left) in order to cope with the Galerkin
approximation of the internal-energy equation in the above proof, which gives
vanishing heat capacity for $\theta<0$ without affecting the sensibility of
the limit problem where one can prove that temperature will actually not go
below 0. The actual internal energy
$\ENG(\FF,\theta)=\psi(\FF,\theta)+\theta\eta(\FF,\theta)$ is then
\begin{align}\label{exa-internal-engr}
\ENG(\FF,\theta)=\dfrac{K_\text{\sc e}}{2\det\FF}\big(\!\det\FF{-}1\big)^2
+\dfrac{G_\text{\sc e}{\rm tr}(\FF\FF^\top)}{(\det\FF)^{1+2/d}}
+\frac{c_{\text{\sc v}}\theta}{(1{+}\ALPHA)\det\!\FF}\Big(\frac{\theta^+}{\theta_0}\Big)^\ALPHA
\ \text{ for }\ \det\FF>0
\end{align}
while $\ENG(\FF,\theta)=+\infty$ for $\det\FF\le0$. In contrast to the referential
internal energy $\Eref$, this actual internal energy contains the term
$K_\text{\sc e}/(2\det\FF)$ and thus exhibits the physically relevant blow-up
\begin{align}\label{exa-internal-engr-blow-up}
\ENG(\FF,\theta)\to+\infty\ \ \ \text{ if }\ \ \det\FF\to0+\,.
\end{align}
This ensures $\det\FF>0$ and thus to prevent local non-interpenetration
directly from the stored energy (uniformly in time) and not only indirectly
by the estimate \eq{est+Fes} which degenerates if the time horizon $T\to\infty$. 
It is also important that the internal energy \eq{exa-internal-engr} is
non-negative even for negative values of temperature. The Cauchy stress $\TT$ is
\begin{align}\label{neo-Hookean-stress}
\TT=\frac{\uppsi_\FF'(\FF,\theta)\FF^\top\!\!\!}{\det\FF}
=\psi_\FF'(\FF,\theta)\FF^\top\!{+}\psi(\FF,\theta)\bbI
=K_\text{\sc e}^{}(\det\FF{-}1)\bbI+2G_\text{\sc e}^{}
\frac{\!\FF\FF^\top\!\!-{\rm tr}(\FF\FF^\top)\bbI/d}{(\det\FF)^{1+2/d}\!}\,;
\end{align}
here we used Cramer's rule $F^{-1}={\rm Cof}\,F^\top\!/\!\det F$ together
with the calculus $\det'(F)={\rm Cof}\,F$. This stress complies trivially with
the energy-controlled-stress condition \eq{ass-stress-control} and $\TT_1\equiv\bm0$.
\end{example}

\begin{remark}[{\sl The energy-controlled Kirchhoff stress}]
\label{rem-control}\upshape
Since additive constants in $\uppsi$ and then in $\ENG$ are actually irrelevant in
the thermodynamical system \eq{Euler-thermo-finite}, the energy control
\eq{ass-stress-control} can equally be written simply as
$|\TT(\FF,\theta)|\le C\ENG(\FF,\theta)$
when assuming a suitable choice of an additive constant in $\uppsi$.
In terms of the referential free energy, this condition means
$|\uppsi_\FF'(\FF,\theta)\FF^\top\!/\!\det\FF|\le C\Eref(\FF,\theta)/\!\det\FF$
with the referential internal energy
$\Eref(\FF,\theta)=\uppsi(\FF,\theta)-\theta\uppsi_\theta'(\FF,\theta)$.
Avoiding a suitable choice of the mentioned additive constants, this can
be written (with a certain tolerance to physical dimensions) as
\begin{align}\label{control-Kirchhoff}
\big|\!\!\!\lineunder{\uppsi_\FF'(\FF,\theta)\FF^\top}{Kirchhoff stress}\!\!\!\big|\le C\big(1+\Eref(\FF,\theta)\big)\,,
\end{align}
as formulated (for the isothermal case) in by J.M.\,Ball
\cite{Ball84MELE,Ball02SOPE}, being devised for usage in the referential
Lagrangian frame and employed e.g.\ in \cite{FraMie06ERCR,MiRoSa18GERV}.
It was shown in \cite[Prop.\,2.3]{Ball02SOPE} that
\eq{control-Kirchhoff} is implied by such energy control of the
Mandel stress $\FF^\top\uppsi_\FF'(\FF,\theta)$. In \cite{MiRoSa18GERV},
such conditions have been used under the name {\it multiplicative stress control}. 
\end{remark}

\begin{remark}[{\sl The free energy blowup \eq{exa-internal-engr-blow-up}}]
\label{rem-neo+}\upshape
Sometimes, \eq{free-neo-Hookean} is expanded by a term $-K_0{\rm ln}(\det\FF)$
with some (presumably small) modulus $K_0>0$ in J/m$^3$=Pa to articulate
the blow-up \eq{exa-internal-engr-blow-up} not only in the actual Eulerian frame
but even in the reference frame. Specifically,
\begin{align}\label{free-neo-Hookean+}
\uppsi(\FF,\theta)=\frac12K_\text{\sc e}^{}\Big(\!\det\FF{-}1{+}\frac{K_0}{K_\text{\sc e}^{}}\Big)^2\!+
G_\text{\sc e}^{}\frac{{\rm tr}(\FF\FF^\top)}{(\det\FF)^{2/d}}
-\frac{c_{\text{\sc v}}\theta}{\!\ALPHA(1{+}\ALPHA)}\Big(\frac{\theta^+}{\theta_0}\Big)^\ALPHA\!
-K_0{\rm ln}(\det\FF)
\end{align}
for $\det\FF>0$ while $\uppsi(\FF,\theta)=+\infty$ for $\det\FF\le0$. Note that
$\uppsi(\cdot,\theta)$ from \eq{free-neo-Hookean+} is minimized on the orbit
${\rm SO}(d)$ as also \eq{free-neo-Hookean} is. This additional term
$-K_0{\rm ln}(\det\FF)$ in \eq{free-neo-Hookean+} also expands the internal
energy \eq{exa-internal-engr} and yields the additional
hydrostatic-pressure-type stress $- K_0\bbI/\!\det\FF$, which complies
with the energy-controlled-stress condition
\eq{ass-stress-control}. Sometimes, a similar term $K_0/\det\FF^m$ with $m>0$ 
in $\uppsi$ is considered instead, referring to the Ogden hyperelastic model and
leading to similar effects and again complying with \eq{ass-stress-control},
cf.\ \cite[Sec.\,5]{FraMie06ERCR}. This is often used in the
Lagrangian formulation. Yet, then the  local non-interpenetration is
ensured differently than here. S{}pecifically, a higher-order gradient
 is considered  in the
stored energy  instead of the dissipative part  and a faster blowup
than only $1/\det\FF$ in Example~\ref{exa-neo} is needed
to exploit the results by T.J.\,Healey and S.\,Kr\"omer \cite{HeaKro09IWSS},
cf.\ also \cite{KruRou19MMCM}.
\end{remark}

\begin{example}[{\sl Volumetric thermal expansion}]\label{exa-neo-exp}
\upshape
Augmenting the free energy \eq{free-neo-Hookean} by the term such
as $-\beta K_\text{\sc e}^{}\theta^+{\rm ln}(\det\FF)$ with $\beta>0$ a
volumetric-expansion coefficient in K$^{-1}$ does not
violate concavity of $\uppsi(\FF,\cdot)$ and allows for modelling
thermal expansion when one realizes that the minimum with respect to
$\det\FF$ of $\uppsi(\cdot,\theta)$ is at $\det\FF=1+\beta\theta+o(\theta)$.
Yet, such  $\uppsi(\cdot,\theta)$ is nonsmooth at $\theta=0$, so we should rather
consider some smooth ansatz, say
\begin{align}\label{free-neo-Hookean-expansion}
\uppsi(\FF,\theta)=\upphi(\FF)
-\frac{c_{\text{\sc v}}\theta}{\!\ALPHA(1{+}\ALPHA)}\Big(\frac{\theta^+}{\theta_0}\Big)^\ALPHA\!
-\beta K_\text{\sc e}^{}\frac{(\theta^+)^2}{\theta^+{+}\theta_{\text{\sc{r}}}}
{\rm ln}(\det\FF)
\end{align}
with the stored energy $\upphi$ from \eq{free-neo-Hookean} and with some fixed
(presumably small) parameter $\theta_{\text{\sc{r}}}>0$.
The last term in \eq{free-neo-Hookean-expansion} contributes to the
referential entropy by $\beta K_\text{\sc e}^{}\theta^+(\theta{+}2\theta_{\text{\sc{r}}}){\rm ln}(\det\FF)/(\theta^+{+}\theta_{\text{\sc{r}}})^2$ so that the 3rd law of
thermodynamics is not corrupted. Its contribution to the referential internal
energy is
$-\beta K_\text{\sc e}^{}(\theta^+)^2(\theta{+}2\theta_{\text{\sc{r}}}){\rm ln}(\det\FF)/(\theta^+{+}\theta_{\text{\sc{r}}})^2$, which is small if
$\theta_{\text{\sc{r}}}>0$ is small and which does not violate non-negativity of
the heat part of the internal energy $\ENT(\FF,\cdot)$. The contribution to the
referential heat capacity 
$2\beta K_\text{\sc e}^{}\theta^+\theta_{\text{\sc{r}}}{\rm ln}(\det\FF)/(\theta^+{+}\theta_{\text{\sc{r}}})^2$ is bounded (uniformly for $\theta_{\text{\sc{r}}}>0$
and, if $\theta_{\text{\sc{r}}}>0$ is small, manifests itself rather only on a small
neighbourhood of the zero temperature. Eventually, its contribution to the Cauchy
stress would be of the hydrostatic-pressure character, namely
$\TT_1=-\beta K_\text{\sc e}^{}(\theta^+\theta/(\theta^+{+}\theta_{\text{\sc{r}}}))\bbI$.
\end{example}

\begin{example}[{\sl A bounded heat capacity}]\label{exa-c-bounded}\upshape
An interesting modification of a concave referential free energy \eq{free-neo-Hookean}
satisfying $\uppsi_\theta'(\FF,0)=0$ and complying with the assumptions
(\ref{Euler-ass}c--f) for $\ALPHA=0$ is
\begin{align}\label{free-neo-Hookean-mod}
\uppsi(\FF,\theta)=\upphi(\FF)
-c_{\text{\sc v}}\Big((\theta^+{+}\theta_{\text{\sc{r}}}){\rm ln}(\theta^+{+}\theta_{\text{\sc{r}}})
-\theta^+\big(1{+}{\rm ln}(\theta_{\text{\sc{r}}})\big)-\theta_{\text{\sc{r}}}{\rm ln}(\theta_{\text{\sc{r}}})\Big)\,,
\end{align}
where $c_{\text{\sc v}}>0$ represents a finite heat capacity asymptotical for
$\theta\to\infty$ for undeformed medium. Like in Example~\ref{exa-neo-exp},
$\theta_{\text{\sc{r}}}>0$ is a fixed (presumably small) temperature-like parameter.
The corresponding actual entropy is
$$
\eta(\FF,\theta)=-\psi_{\theta}'(\FF,\theta)
=\frac{c_{\text{\sc v}}}{\det\FF}\,{\rm ln}\Big(1+\frac{\theta^+}{\theta_{\text{\sc{r}}}}\Big)
$$
and the actual internal energy 
$\ENG(\FF,\theta)=\psi(\FF,\theta)+\theta\eta(\FF,\theta)$ is
\begin{align}\label{internal-neo-Hookean-mod}
\ENG(\FF,\theta)
=\frac{\upphi(\FF)}{\det\FF}+
\frac{c_{\text{\sc v}}}{\det\FF}\bigg(\theta^++\theta_{\text{\sc{r}}}{\rm ln}\Big(1+\frac{\theta^+}{\theta_{\text{\sc{r}}}}\Big)\bigg)\,.
\end{align}
Realizing that 
$\uppsi_{\theta\theta}''(\FF,\theta)=-(c_{\text{\sc v}}/\theta_{\text{\sc{r}}})/(1{+}\theta^+/\theta_{\text{\sc{r}}})$,
we obtain a bounded increasing actual heat capacity $c=c(\FF,\cdot)$ with
$c(\FF,0)=0$ and $c(\FF,+\infty)=c_{\text{\sc v}}/\!\det\FF$, namely
\begin{align}
c(\FF,\theta)=\frac{c_{\text{\sc v}}\theta^+}{(\theta^+{+}\theta_{\text{\sc{r}}})\det\FF}\,.
\end{align}
A comparison of this model with Example~\ref{exa-neo} is in Figure~\ref{fig-free-energy}.
\begin{figure}[ht]\hspace*{-.2em}
\begin{tikzpicture}  
    \begin{axis}[width=6.5 cm,height=6.4 cm,
        xmin=-1.2,xmax=4.5,ymin=-42,ymax=2,  
        clip=true,  
        axis lines=center,  
        grid = major,  
        ytick={-40, -30, -20, -10, 0},  
        xtick={-1, 0, 1, 2, 3, 4},  
        xlabel=$\theta$,
        ylabel={$\uppsi=\uppsi(\theta)$\hspace*{-1em}},  
        every axis y label/.style={at=(current axis.above origin),anchor=south},  
        every axis x label/.style={at=(current axis.right of origin),anchor=west},  
      ]  
\addplot[very thick,domain=-1:5,samples=30]{0.2*ln(0.2)-(max(0,x)+0.2)*ln(max(0,x)+0.2)+max(0,x)*(1+ln(0.2))};
\addplot[very thick,domain=-1:5.5,samples=30,dotted]{-1*max(0,x)^(1+0.2)/(0.2*(1+0.2))};
\addplot[very thick,domain=-1:5.5,samples=30,dashed]{-1*max(0,x)^(1+0.05)/(0.05*(1+0.05))};
  \end{axis}
  \end{tikzpicture}
\begin{tikzpicture}  
    \begin{axis}[width=6.5 cm,height=6.3 cm,
        xmin=-1.3,xmax=5.5,ymin=-.5,ymax=5.2,  
        clip=true,  
        axis lines=center,  
        grid = major,  
        ytick={0, 1, 2, 3, 4, 5},  
        xtick={-1, 0, 1, 2, 3, 4, 5},  
        xlabel=$\theta$,
        ylabel={$\Eref=\Eref(\theta)$\hspace*{-1em}},  
        every axis y label/.style={at=(current axis.above origin),anchor=south},  
        every axis x label/.style={at=(current axis.right of origin),anchor=west},  
      ]  
\addplot[very thick,domain=-1:5.5,samples=70]{max(0,x)+0.2*ln(max(0,x)/0.2+1))};
\addplot[very thick,domain=-1:5.5,samples=30,dotted]{1*max(0,x)^(1+0.2)/(1+0.2)};
\addplot[very thick,domain=-1:5.5,samples=30,dashed]{1*max(0,x)^(1+0.05)/(1+0.05)};
  \end{axis}
  \end{tikzpicture}
\begin{tikzpicture}
    \begin{axis}[width=6.5 cm,height=6.15 cm,
        xmin=-1.3,xmax=4.6,ymin=-.1,ymax=1.35,  
        clip=true,  
        axis lines=center,  
        grid = major,  
        ytick={0, 0.2, 0.4, 0.6, 0.8, 1, 1.2},  
        xtick={-1, 0, 1, 2, 3, 4},  
        xlabel=$\theta$,
        ylabel={${\rm c}={\rm c}(\theta)$\hspace*{-1em}},  
        every axis y label/.style={at=(current axis.above origin),anchor=south},  
        every axis x label/.style={at=(current axis.right of origin),anchor=west},  
      ]  
\addplot[very thick,domain=-1.1:5,samples=70]{max(0,x)/(max(0,x)+0.2)};
\addplot[very thick,domain=-1.1:5,samples=70,dotted]{1*max(0,x)^0.2};
\addplot[very thick,domain=-1.1:5,samples=70,dashed]{1*max(0,x)^0.05};
\node at (axis cs:1.,.08) {$\larger{\larger\boldsymbol\downarrow}$};
\node at (axis cs:1.,.92) {$\larger{\larger\boldsymbol\uparrow}$};
\node at (axis cs:1.09,0.19) {$\boldsymbol{\theta_0}$};  
            \end{axis}  
  \end{tikzpicture}\hspace*{-.0em}
\hspace*{-.5em}
\caption{{\sl A comparison of Examples~\ref{exa-neo} and \ref{exa-c-bounded}
with $c_{\text{\sc v}}=1$ ignoring the mechanical part, i.e.\ $\upphi=0$. 
The former example is depicted for $\ALPHA=0.05$ (dashed lines) and for
$\ALPHA=0.2$ (dotted lines) while the latter  example is depicted for
$\theta_{\text{\sc{r}}}=0.2$ (solid lines). Both examples yield a similar internal
energy $\ENG$ and a similar heat capacity $c$ for
$\ALPHA\to0$ and $\theta_{\text{\sc{r}}}\to0$. }}
\label{fig-free-energy}\end{figure}
\end{example}

\begin{example}[{\sl The phase transition in shape-memory alloys}]\label{exa-SMA}
\upshape
The energy $\uppsi$ from \eq{free-neo-Hookean} has a single well in the sense
that the minimum $\uppsi(\cdot,\theta)$ is attained on the orbit
${\rm SO}(d)$. A multi-well modification is used to model a so-called
martensitic transformation in shape-memory alloys on the level of single
crystals, cf.\ \cite{Roub04MMES} for an overview.
One considers, beside one so-called austenitic well as in \eq{free-neo-Hookean},
also $L$ lower-symmetrical wells to models variants of so-called martensite.
In particular, $L=3$ for the tetragonal martensite while $L=6$ in case of
 the orthorhombic martensite or $L=12$ for the monoclinic martensite.
The stress-free geometric configuration of particular phases is characterized
by the distortion $\mathbb F_\ell$ with $\ell=0,1,...,L$. 
For particular phases, the ansatz \eq{free-neo-Hookean} is the modified as
\begin{align}\nonumber
\uppsi_\ell(\FF,\theta)=\upphi_\ell(\FF)
-\frac{{\rm c}_\ell\theta}{\!\ALPHA(1{+}\ALPHA)}\Big(\frac{\theta^+}{\theta_0}\Big)^\ALPHA
\ \ \text{ with }\ \ \upphi_\ell(\FF)=\frac12K_\ell^{}\big(\det\FF{-}1\big)^2\!+
G_\ell^{}
\frac{{\rm tr}(\mathbb F_\ell^{-\top}\!\FF\FF^\top\!\mathbb F_\ell^{-1})\!}{(\det\FF)^{2/d}}\,.
\end{align}
The minimum of $\uppsi_\ell(\cdot,\theta)$ is attained on the orbit ${\rm SO}(d)\FF_\ell$. For austenite, $\mathbb F_0=\bbI$ while $\mathbb F_\ell\ne\bbI$ for $\ell=1,...,L$ and, since martensitic transformation is essentially isochoric, 
$\det\mathbb F_\ell=1$ for $\ell=1,...,L$.
The overall referential  multi-well free energy is 
\begin{align}\label{SMA}
\uppsi(\FF,\theta)=-\varkappa\,{\rm ln}\bigg(\sum_{\ell=0}^L
{\rm e}^{\displaystyle{-\uppsi_\ell(\FF,\theta)/\varkappa}}\bigg)\,,
\end{align}
where $\varkappa$ is a constant with the physical dimension J/m$^3$=Pa; one can think about
the Boltzmann constant (related per unit volume) times a reference temperature.
Such  $\uppsi(\cdot,\theta)$ exhibits the {\it multi-well} character and is backed up by
statistical physics. The so-called shape-memory effect is modelled by considering
${\rm c}_0>{\rm c}_1=...={\rm c}_L$, which causes that the well of austenite falls
energetically faster
than the well of martensitic variants within increasing temperature, which makes
austenite energetically dominant at high temperatures while martensite is dominant
at lower temperatures. The Cauchy stress  corresponding to \eq{SMA} is
\begin{align}\label{SMA-stress}
\TT(\FF,\theta)=\dfrac{\sum_{\ell=0}^L{\rm e}^{-{\uppsi_\ell(\FF,\theta)/\varkappa}}
\upphi_\ell'(\FF)}{\sum_{\ell=0}^L{\rm e}^{-{\uppsi_\ell(\FF,\theta)/\varkappa}}\det\FF}\FF^\top
\end{align}
while the actual entropy 
$\eta={-}\psi_\theta'(\FF,\theta)={-}\uppsi_\theta'(\FF,\theta)/\det\FF$ can be
evaluated as
\begin{align}\label{SMA-entropy}
\eta(\FF,\theta)=\dfrac{\sum_{\ell=0}^L
{\rm e}^{-\psi_\ell(\FF,\theta)/\varkappa}\eta_\ell(\FF,\theta)}{\sum_{\ell=0}^L
{\rm e}^{-\psi_\ell(\FF,\theta)/\varkappa}}\ \ \ \text{ with }\ \ \eta_\ell(\FF,\theta)=
-\frac{{\rm c}_\ell}{\ALPHA\det\FF}\Big(\frac{\theta^+}{\theta_0}\Big)^\ALPHA\,,
\end{align}
and the actual internal energy $\Eng=\psi+\theta\eta$ is
\begin{align}\label{SMA-internal}
\Eng=\ENG(\FF,\theta)=\dfrac{\sum_{\ell=0}^L
{\rm e}^{-\uppsi_\ell(\FF,\theta)/\varkappa}\,\ENG_\ell(\FF,\theta)}{\sum_{\ell=0}^L
{\rm e}^{-\uppsi_\ell(\FF,\theta)/\varkappa}}\ \ \text{ with }\
\ENG_\ell(\FF,\theta)=\psi_\ell(\FF,\theta)+\theta\eta_\ell(\FF,\theta)\,,
\end{align}
where naturally $\psi_\ell(\FF,\theta)=\uppsi_\ell(\FF,\theta)/\!\det\FF$.
Let us remark that the above isotropic neo-Hookean $\upphi_\ell$  
is typically considered rather anisotropic, reflecting lower-symmetrical
character of particular martensitic variants for $\ell=1,...,L$.
\end{example}

\begin{remark}[{\sl Nonphysical temperature-independent heat capacity}]
\label{rem-nonphysical}\upshape
Often, the specific heat capacity is considered nearly constant
for temperatures far from the absolute zero (so-called
{\it Dulong–Petit}'s {\it law}), which suggests $\ALPHA>0$ in \eq{free-neo-Hookean}
or $\theta_{\text{\sc{r}}}>0$ in \eq{free-neo-Hookean-expansion} or in
\eq{free-neo-Hookean-mod} to be very small, cf.\ also
Figure~\ref{fig-free-energy}-right. On the other hand, the above examples do not
work for $\ALPHA=0$ nor for $\theta_{\text{\sc{r}}}=0$.
To obtain the truly temperature-independent heat capacity, the thermal term in
the free energy \eq{free-neo-Hookean} or in \eq{free-neo-Hookean-mod} is to be
modified as
\begin{align}\label{non-physical}
\uppsi(\FF,\theta)=\upphi(\FF)-c_{\text{\sc v}}\theta{\rm ln}\Big(\frac\theta{\theta_0}\Big)\,.
\end{align}
Although this is the frequently considered ansatz in engineering and
mathematical literature, the corresponding actual entropy
$c_{\text{\sc v}}{\rm ln}(\theta/\theta_0)/\det\FF$
obviously conflicts with the 3rd law of thermodynamics.
\end{remark}


\end{sloppypar}
\end{document}